\documentclass[11pt]{amsart}
\usepackage{xypic}
\usepackage{amscd}
\usepackage{epsf}
\xyoption{ps}\xyoption{dvips} 
\newbox\mybox 
\def\overtag#1#2#3{\setbox\mybox\hbox{$#1$}\hbox to
  0pt{\vbox to 0pt{\vglue-#3\vglue-\ht\mybox\hbox to \wd\mybox
      {\hss$\ss#2$\hss}\vss}\hss}\box\mybox}
\def\undertag#1#2#3{\setbox\mybox\hbox{$#1$}\hbox to 0pt{\vbox to
    0pt{\vglue#3\vglue\ht\mybox\hbox to \wd\mybox
      {\hss$\ss#2$\hss}\vss}\hss}\box\mybox} 
\def\lefttag#1#2#3{\hbox to 0pt{\vbox to 0pt{\vss\hbox to
      0pt{\hss$\ss#2$\hskip#3}\vss}}#1} 
\def\righttag#1#2#3{\hbox to 0pt{\vbox to 0pt{\vss\hbox to
      0pt{\hskip#3$\ss#2$\hss}\vss}}#1} 
\let\ss\scriptstyle

\def\Dot{\lower.2pc\hbox to 2pt{\hss$\bullet$\hss}}
\def\Circ{\lower.2pc\hbox to 2pt{\hss$\circ$\hss}}
\def\Vdots{\raise5pt\hbox{$\vdots$}}
\newcommand\lineto{\ar@{-}}
\newcommand\dashto{\ar@{--}}
\newcommand\dotto{\ar@{.}}
\def\WDN#1.{\marginpar{\rightskip=0pt plus2\hsize#1}} 
\def\mt#1_#2{#1\times_{#2}S^1}
\def\pmt#1_#2{(#1\times_{#2}S^1)}
\newcommand{\interior}[1]{\vphantom{\vbox{\vbox
 {\hbox{$\scriptstyle\circ$}\vskip.3pt}\nointerlineskip\hbox{$#1$}}}
 \vbox{\vbox to 0pt{\vss\hbox{\hskip3pt$\scriptstyle\circ$}
 \vskip.4pt}\nointerlineskip\hbox{$#1$}}}

\newcommand\Interior{\operatorname{int}}
\newcommand\var{\operatorname{var}}
\newcommand\Ker{\operatorname{Ker}}

\renewcommand\Im{\operatorname{Im}}
\newcommand\fbar{\overline f} 

\newcommand\C{{\mathbb C}}
\newcommand\Z{{\mathbb Z}}
\newtheorem{theorem}{Theorem}[section]
\newtheorem{lemma}[theorem]{Lemma}
\newtheorem*{scholium}{Scholium}
\newtheorem{proposition}[theorem]{Proposition}
\newtheorem{corollary}[theorem]{Corollary}

\theoremstyle{definition}
\newtheorem{definition}[theorem]{Definition}

\newtheorem*{terminology}{Terminology}

\begin{document}
\title[Unfolding polynomial maps at infinity] {Unfolding polynomial
  maps at infinity}
\author{Walter D. Neumann} 
\address{Department of Mathematics and Statistics\\The University of
  Melbourne\\Vic 3010, Australia}
\email{neumann@ms.unimelb.edu.au} 
\author{Paul Norbury}
\address{Department of Mathematics and Statistics\\The University of
  Melbourne\\ Vic 3010, Australia}
\email{norbs@ms.unimelb.edu.au} \keywords{}  
\thanks{This research was supported by the Australian  Research Council}
\maketitle

\section{Introduction}
\label{sec:intro}

Let $f\colon \C^n\to\C$ be a polynomial map.  The polynomial describes
a family of complex affine hypersurfaces $f^{-1}(c)$, $c\in \C$.  The
family is locally trivial, so the hypersurfaces have constant
topology, except at finitely many \emph{irregular} fibers $f^{-1}(c)$
whose topology may differ from the generic or \emph{regular} fiber of
$f$.  

We would like to give a full description of the topology of this
family in terms of easily computable data. This paper describes some
progress.

We will restrict mostly to the case that $f$ has only isolated
singularities. We show that the data then needed are local
monodromy maps obtained by transporting a fixed generic fiber $F$
around the irregular fibers, and the ``Milnor fibers'' of the singular
points and singularities at infinity of $f$: these are certain
submanifolds of $F$ that describe the loss of topology at irregular
fibers.

It is convenient to subdivide this necessary data as follows. For each 
irregular fiber we need
\begin{itemize}
\item the Milnor fibers associated with it;
\item the local monodromy for the irregular fiber restricted to each
  Milnor fiber;
\item the embeddings of the Milnor fibers into a fixed ``reference''
  regular fiber $F$.
\end{itemize}

The first two items are local ingredients, while the third is global.
We are able to give complete computation of the local ingredients for
$n=2$ (see Theorem \ref{th:dim2}).  The remaining problem for $n=2$ is
therefore the third item, although we obtain enough constraints that
the complete topology can be sometimes be deduced.

We use the Brian\c con polynomial as an illustrative example. In this
case our general results quickly yield the previous homological
monodromy computations of Artal-Bartolo, Cassou-Nogues, and Dimca
\cite{artal-cassou-dimca} and Dimca and Nemethi \cite{dimca-nemethi}.
Since our computations are geometric, we obtain sharper information
(action on the intersection form, etc.).  In fact, our computations
strongly suggest a candidate description for the complete topology.
However, this description remains conjectural, so the example
illustrates both the strengths and current limitations of our
approach.

As in the paper \cite{neumann-rudolph}, we also give a detailed
description on the level of homology. That paper shows that, if $f$ is
\emph{good}, that is, it has only isolated singularities and no
singularities at infinity, then for $n>3$ the homological data
describes the complete topology.  Our results prepare the ground for
an extension of this result to the case that there are singularities
at infinity (at least if the singularities at infinity are also
isolated), by extending the topological model of
\cite{neumann-rudolph} to allow singularities at infinity.

Before we describe our results in more detail, we say a bit more about
the basic set-up.

A fiber can be irregular in two ways:
\begin{itemize}
\item it may have singularities;
\item it may be ``irregular at infinity'': failure of local triviality
  near the given fiber occurs outside arbitrarily large compact sets.
\end{itemize}
Of course a fiber can be both singular and irregular at infinity; this
will, for instance, always be the case if the fiber has non-isolated
singularities.  There have been many papers dealing with algebraic
conditions that imply regularity at infinity (M-tameness of
\cite{nemethi-zaharia} which asks that the given fiber and all nearby
fibers be transverse to sufficiently large spheres; $\rho$-regularity
of \cite{tibar} and \cite{tibar-Wallpaper} which generalizes this to
allow non-round spheres; the stronger t-regularity of
\cite{siersma-tibar}, equivalent to the ``Malgrange condition'' of
\cite{pham}).  For $n=2$ these conditions are all equivalent to
regularity at infinity (see \cite{durfee}, \cite{tibar-Wallpaper}),
but in higher dimensions they are not mutually equivalent. For our
purposes the weakest concept of irregularity at infinity, namely the
topological one given above, suffices.

In \cite{neumann-norbury} we described results that hold under no
conditions on singularities. As already mentioned, we will here be
particularly interested in fibers which have only isolated
singularities, but no restriction on singularities at infinity.  There
have been several recent studies of polynomial maps for which the
singularities at infinity are also restricted to be isolated in an
appropriate sense (e.g., \cite{artal-cassou-dimca}, \cite{broughton1},
\cite{broughton2}, \cite{siersma-tibar}, \cite{parusinski},
\cite{parusinski2}).  However most of our results are new even for
this case.  In particular, the second part of this paper deals only
with the case of dimension $n=2$, in which case singularities at
infinity are always isolated.

At an irregular fiber there is a loss of topology compared with the
regular fiber.  For an isolated singularity the change in topology is
captured by the Milnor fiber \cite{milnor} of the singularity.  If $f$
is ``good'' (no singularities at infinity) then these Milnor fibers
account for all the topology of the regular fiber in a way that is
made precise by Neumann and Rudolph \cite{neumann-rudolph}.

For a singularity at infinity there is a ``Milnor fiber at infinity,''
first described by Suzuki \cite{suzuki} for $n=2$.  One point of this
paper is to extend the theory of \cite{neumann-rudolph} to encompass
the Milnor fibers at infinity also, and to show that the new
ingredients --- the topology of the singularities at infinity of an
irregular fiber, as encoded by the Milnor fibers at infinity and the
monodromy maps on these Milnor fibers --- are recoverable for $n=2$ in
full from the link at infinity of the irregular fiber.


\section{Homological results}\label{sec:homological}
To clarify what we mean by ``loss of topology at an irregular fiber''
we start with a simple homological result which is true under no
assumption on $f$.

\begin{definition}\label{def:vanishing}
  For any fiber $f^{-1}(c)$ of $f\colon\C^n\to\C$ choose $\epsilon$
  sufficiently small that all fibers $f^{-1}(c')$ with $c'\in
  D^2_\epsilon(c)-\{c\}$ are regular ($D^2_\epsilon(c)$ is the closed
  disk of radius $\epsilon$ about $c$) and let
  $N(c):=f^{-1}(D^2_\epsilon(c))$. Let $F=f^{-1}(c')$ be a regular
  fiber in $N(c)$. Then
  \begin{align*}
    V_q(c)&:= \Ker\bigl(H_q(F;\Z)\to H_q(N(c);\Z)\bigr)
  \end{align*}
  is the group of \emph{vanishing $q$-cycles} for $f^{-1}(c)$.
\end{definition}

Let $\Sigma$ be the set of irregular values of $f$. Choose a regular
value $c_0$ for $f$ and paths $\gamma_c$ from $c_0$ to $c$ for each
$c\in\Sigma$ which are disjoint except at $c_0$.  We use these paths
to refer homology of a regular fiber near one of the irregular fibers
$f^{-1}(c)$ to the homology of the ``reference'' regular fiber
$F=f^{-1}(c_0)$.  From now on $F$ will always mean this particular
regular fiber.

\begin{theorem}[\cite{broughton1}, \cite{neumann-norbury}]\label{th:van}
  For $q>0$ the maps $V_q(c)\to H_q(F;\Z)$ induce an isomorphism
  \begin{equation*}
\bigoplus_{c\in\Sigma} V_q(c) \cong H_q(F;\Z) .
  \end{equation*}
  Moreover, the map $H_q(F;\Z)\to H_q\bigl(N(c);\Z\bigr)$ with kernel
  $V_q(c)$ is surjective, so
  $H_q\bigl(N(c);\Z\bigr)\cong\bigoplus_{c'\in\Sigma-\{c\}}V_q(c')$.\qed
\end{theorem}
\noindent 
Thus the group of vanishing cycles measures homologically the ``loss
of topology'' at an irregular fiber, and these groups account for all
the homology of the regular fiber.

We now restrict to a fiber $f^{-1}(c)$ with at most isolated
singularities (but possibly singular at infinity). The ``nonsingular
core'' of $f^{-1}(c)$ is obtained by intersecting $f^{-1}(c)$ with a
very large ball and then removing small regular neighborhoods of its
singularities.  More precisely,
there is a radius $R(c)$ such that for any $r\ge R(c)$ the sphere
$S^{2n-1}_r\subset\C^n$ of radius $r$ about the origin intersects
$f^{-1}(c)$ transversally.  Choose any $r\ge R(c)$ and denote
$F^{co}(c):= f^{-1}(c)\cap D^{2n}_r(0)$, where $D^{2n}_r(0)$ is the
disk of radius $r$ about the origin in $\C^n$.  This $F^{co}(c)$ is
the \emph{compact core} of the fiber $f^{-1}(c)$.  The fiber
$f^{-1}(c)$ is topologically the result of adding an open collar to
the boundary of $F^{co}(c)$.
Any singularities of $f$ on
$f^{-1}(c)$ lie on $F^{co}(c)$ and we remove small open regular
neighborhoods of these singularities to form a boundaried
$2n$-manifold $F^{ns}(c)$, the \emph{non-singular core} of
$f^{-1}(c)$.

There may be boundary components of $F^{co}(c)$ outside of which the
topology of nearby fibers of $f$ is the same as that of $f^{-1}(c)$,
in the sense that $f$ restricted to the appropriate component of
$f^{-1}\bigl(D^2_\delta(c)\bigr)\cap\bigl(\C^n-\interior
D^{2n}_r(0)\bigr)$ gives a locally trivial fibration for $r\ge R(c)$
and $\delta$ sufficiently small.  We call such boundary components
\emph{regular} and call the other boundary components of $F^{co}(c)$
\emph{irregular}. Thus $f^{-1}(c)$ is regular at infinity if and only
if all boundary components of $F^{co}(c)$ are regular.

By standard arguments (see Section \ref{sec:alldim}) we can embed
$$\phi\colon F^{ns}(c)\times D^2_\epsilon(c)\hookrightarrow \C^n$$
(if
$\epsilon$ is small enough) so that $\phi(x,c)=x$ for $x\in F^{ns}(c)$
and $f\circ\phi$ is the projection to $D^2_\epsilon(c)$.  By
restricting $\phi$ to $F^{ns}(c)\times\{c'\}$ we thus get an embedding
of $F^{ns}(c)$ into a nearby regular fiber $F'=f^{-1}(c')$ of $f$.  The
complement $F'-\interior F^{ns}(c)$ then consists of the
disjoint union of the Milnor fibers of the singularities of $f$ on
$f^{-1}(c)$ and certain non-compact pieces.  These non-compact pieces
will be half-open collars on the regular boundary components of
$F^{co}(c)$ and other pieces which meet $F^{ns}(c)$ at irregular
boundary components of $F^{co}(c)$. We call the latter the
\emph{Milnor fibers at infinity} for $f^{-1}(c)$.

\begin{terminology}
  By \emph{Milnor fibers of $f$} we will mean all Milnor fibers of
  isolated singularities of $f$ and all Milnor fibers at infinity.  If
  we want to emphasize that a Milnor fiber is not at infinity we will
  call it a \emph{finite Milnor fiber}.  Let $F_1(c),\dots,F_{s_c}(c)$
  and $F_{s_c+1}(c),\dots,F_{t_c}(c)$ be all the Milnor fibers at
  infinity respectively finite Milnor fibers for $f^{-1}(c)$.
\end{terminology}

We will consider the Milnor fibers to lie in our standard regular
fiber $F$, by transporting the fiber $F'=f^{-1}(c')$ along the path
$\gamma_c$.  Topologically, the fiber $f^{-1}(c)$ results from $F$ by
collapsing each finite Milnor fiber to a point and removing each
Milnor fiber at infinity.  This is the sense in which the Milnor
fibers capture the loss of topology of the fiber $f^{-1}(c)$.  We can
use the local monodromy to relate this to vanishing cycles.

For each irregular fiber $f^{-1}(c)$, by transporting a nearby regular
fiber $F'$ in a small loop around the fiber $f^{-1}(c)$ we get a local
monodromy map $F'\to F'$. By using the path $\gamma_c$ to refer this
monodromy to the reference regular fiber $F=f^{-1}(c_0)$ we consider
it as a map $h(c)\colon F\to F$.  This monodromy map is well
defined up to isotopy.  

If $f^{-1}(c)$ has isolated singularities, the local monodromy
$h(c)\colon F\to F$ can be normalized to be the identity on the image
in $F$ of the non-singular core of $f^{-1}(c)$ (use the above
embedding $\phi\colon F^{ns}(c)\times D^2\hookrightarrow \C^n$).
Thus, $h(c)$ restricts to a local monodromy map on each Milnor fiber
$F_i(c)$. We denote this local monodromy map
\begin{equation*}
h(c)\colon F_i(c)\to F_i(c) 
\end{equation*}
also by $h(c)$, or simply $h$. 

This map is the identity on $\partial F_i(c)$, so it induces a map in
homology called the \emph{variation} (introduced by \cite{lamotke},
but with different sign convention)
\begin{equation*}
  \var\colon H_q\bigl(F_i(c),\partial
F_i(c)\bigr)\to H_q\bigl(F_i(c)\bigr),
\end{equation*}
 obtained by taking a relative cycle $C$ to
the closed cycle $C-h_\sharp C$.  Let 
\begin{equation*}
  (v_i)_q\colon
H_q\bigl(F_i(c),\partial F_i(c)\bigr)\to H_q(F)
\end{equation*}
be the composition of variation
with the map $H_q(F_i(c))\to H_q(F)$ induced by inclusion.

\begin{theorem}\label{th:homology}
  For $q\ge1$ the maps $(v_i)_q\colon H_q(F_i(c),\partial F_i(c))\to H_q(F)$
  are injective and induce an isomorphism 
  \begin{equation*}
    \bigoplus_{i=1}^{t_c}
  H_q\bigl(F_i(c),\partial F_i(c)\bigr)\stackrel{\cong
}{\longrightarrow} V_q(F)
  \end{equation*}
to the subgroup $V_q(F)\subset H_q(F)$ of vanishing cycles.
\end{theorem}

We will prove this theorem in section \ref{sec:alldim}. We first
refine it and Theorem \ref{th:van} by describing how monodromy,
intersection form, and Seifert form relate to these sum decompositions
of homology.

\section{Monodromy and Seifert form in homology}\label{sec:monodromy}
\begin{theorem}
  \label{th:monodromy}
  With respect to the sum decomposition in Theorem \ref{th:van} the
  map in homology induced by the local monodromy $h(c)$ has the form
$$
\begin{pmatrix}
  I&0&\dots&0&0&\dots&0\\
\vdots&\vdots&&&&\vdots\\
*&*&\dots&h_c&*&\dots&*\\
\vdots&\vdots&&&&\vdots\\
0&0&\dots&0&0&\dots&I
\end{pmatrix}$$
where $h_c$ is the restriction of $H_*(h(c))$ to $V_q(c)$.

If $f^{-1}(c)$ has only isolated singularities, then $h_c$ respects
the sum decomposition of Theorem \ref{th:homology} and thus
has block form
$$
\begin{pmatrix}
  h_{c,1}&0\dots&0\\
0&{h_c,2}\dots&0\\
\vdots&\vdots&&\vdots\\
0&0&\dots&h_{c,t_c}
\end{pmatrix}
$$
where $h_{c,i}$ is induced by the local monodromy $h(c)\colon
F_i(c)\to F_i(c)$.
\end{theorem}

The \emph{link at infinity} of a fiber $f^{-1}(c)$ with isolated
singularities is the link
$\bigl(S^{2n-1},S^{2n-1}\cap f^{-1}(c)\bigr)$, where $S^{2n-1}$ is any
sphere around the origin of radius greater than the number $R(c)$
mentioned above.  For a regular fiber this link is independent, up to
equivalence, of the choice of fiber and is called the \emph{regular
  link at infinity for $f$}.  A standard construction shows that a
Seifert surface of this regular link at infinity $(S^{2n-1},L_{reg})$
is diffeomorphic to the compact core $F^{co}$ of a regular fiber (for
$n=2$ this gives the minimal Seifert surface, and this minimal Seifert
surface is unique, see \cite{neumann-inv}). 

The Seifert linking form
$H_q(F^{co})\otimes H_{2n-2-q}(F^{co})\to\Z$ on the homology of a
Seifert surface is a useful invariant of a link.

For each $F_i(c)$ we have a form:
\begin{equation}
\begin{split}
  \label{eq:seifert}
  (L_{c,i})_q\colon H_q(F_i(c),\partial F_i(c))\otimes& H_{2n-2-q}(F_i(c),\partial F_i(c))\to \Z,\\
  \alpha\otimes \beta&\mapsto \var(\alpha).\beta,
\end{split}
\end{equation}
where the dot represents intersection form 
\begin{equation*}
  H_q\bigl(F_i(c)\bigr)\otimes
H_{2n-2-q}\bigl(F_i(c),\partial F_i(c)\bigr)\to \Z.
\end{equation*}
\begin{theorem}\label{th:seifert}
Suppose $f$ has only isolated singularities, so Theorems \ref{th:van} and
  \ref{th:homology} give a direct sum decomposition
  \begin{equation*}
    \bigoplus_{c\in\Sigma}\bigoplus_{i=1}^{t_c} 
H_q\bigl(F_i(c),\partial F_i(c)\bigr)\stackrel{\cong}{\longrightarrow}
H_q(F)\cong H_q(F^{co}).
  \end{equation*}
  
  If we order the set $\Sigma$ so that the paths $\gamma_c$, $c\in
  \Sigma$, depart the point $c_0$ in anti-clockwise order, then with
  respect to this direct sum decomposition the Seifert form of the
  regular link at infinity has lower triangular block form, with
  diagonal blocks given%
{} by the forms
  $L_{c,i}$ of {\rm(\ref{eq:seifert})}, and with off-diagonal blocks
  equal to zero for pairs of summands with the same $c$:
$$
\begin{pmatrix}
&&\dots&0&0&\dots&0&0&\dots&0\\
\vdots&\vdots&&&&&&&&\vdots\\
*&*&\dots&L_{c,1}&0&\dots&0&0&\dots&0\\
*&*&\dots&0&L_{c,2}&\dots&0&0&\dots&0\\
*&*&\dots&\vdots&\vdots&&\vdots&0&\dots&0\\
*&*&\dots&0&0&\dots&L_{c,t_c}&0&\dots&0\\
\vdots&\vdots&&&&&&&&\vdots\\
*&*&\dots&*&*&\dots&*&&&\\
\end{pmatrix}$$
\end{theorem}
If $f$ is good, that is, there are no Milnor fibers at infinity, these
results are in \cite{neumann-rudolph} in slightly different
formulation, as follows.  The variation map $\var\colon
H_q\bigl(F_i(c),\partial F_i(c)\bigr)\to H_q(F_i(c))$ is an
isomorphism for the Milnor fiber of an isolated
singularity\footnote{This is a general fact about fibered links, see
  \cite{lamotke}.}. Moreover, in this case the homology $H_q(F_i(c))$
vanishes for $q\ne n-1$ by \cite{milnor}.  Thus, we can replace
$(v_i)_q$ in Theorem \ref{th:homology} by the map
$H_{n-1}\bigl(F_i(c)\bigr)\to H_{n-1}(F)$ induced by inclusion when
$q=n-1$ (and ignore it when $q\ne n-1$). Moreover, in this case the
form of (\ref{eq:seifert}) is the Seifert form pulled back to
$H_*(F_i,\partial F_i)$ via the variation isomorphism. Thus Theorem
\ref{th:homology} can be formulated in terms of the Seifert forms of
the singularities of $f$, which is the form in which these results
were given in \cite{neumann-rudolph}.

In this case that $f$ is good the Seifert form is a particularly
strong invariant (see \cite{neumann-rudolph}): the local homological
monodromies of Theorem \ref{th:monodromy} are all computable
from the the above block decomposition of the Seifert form on
$H_{n-1}(F)$, and for $n>3$ the complete topology of
$f\colon\C^n\to\C$ is determined by
this block decomposition of the Seifert form.

To describe further relations on the above block decompositions, we
suppose the irregular values are numbered $c_1, \dots,c_k$ in the
order they occur in Theorem \ref{th:seifert}. We abbreviate $h_{c_j}$
as $h_j$ and write the decompositions of monodromy and Seifert form of
Theorems \ref{th:monodromy} and \ref{th:seifert} as:
\begin{align}\label{eq:H}
H_*(h(c_j))&=
\begin{pmatrix}
I&0&\dots&0&0&\dots&0\\
\vdots&&&&&&\vdots\\
0&0&\dots&I&0&\dots&0\\
h_{j1}&h_{j2}&\dots&h_{j,j-1}&h_j&\dots&h_{jk}\\
\vdots&&&&&&\vdots\\
0&0&\dots&0&0&\dots&I
\end{pmatrix} \\ \label{eq:L}
L&=
\begin{pmatrix}
L_1&0&\dots&0&0&\dots&0\\
L_{21}&L_2&\dots&0&0&\dots&0\\
\vdots&&&&&&\vdots\\
L_{j1}&L_{j2}&\dots&L_{j}&0&\dots&0\\
\vdots&&&&&&\vdots\\
L_{k1}&L_{k2}&\dots&l_{kj}&\dots&\dots&L_k
\end{pmatrix}~.
\end{align}
As described above, $h_j$ and $L_j$ may decompose further as the block
sums of the $h_{c_j,i}$ respectively $L_{c_j,i}$, $i=1,\dots,t_{c_j}$.

It is a standard result that the intersection form $S$ on $H_*(F)$ may
be written
$$S=L-L^t$$
where $L$ is the Seifert form discussed above, and $L^t$
is the appropriate graded transpose,
$L^t(x,y)=(-1)^{(p+1)(q+1)}L(y,x)$ if $x\in H_q(F)$ and $y\in H_p(F)$
(see, e.g., \cite{durfee-var}; in \cite{neumann-rudolph} this formula
is mistakenly written $S=L+L^t$, and the second instance of
$(-1)^{(s+1)(n-s)}$ on the same page 418 should be $(-1)^{s+1}$). The
intersection form is, of course, preserved by all the local
monodromies $h(c_j)$.  The Seifert form, on the other hand, is only
preserved by the ``monodromy at infinity''
$h(\infty)=h(c_k)h(c_{k-1})\dots h(c_1)$.  There are, nevertheless,
some relationships between local monodromy and Seifert form which can
give useful constraints.  The following generalizes Theorem 3.5 of
\cite{neumann-rudolph}.
\begin{theorem}\label{th:relations}
  With notation as above,
\begin{align*}
  L_jh_j&=L_j^t,  \\
L_ih_{ij}&=L_{ji}^t&&\text{for $i<j$,}\\
L_ih_{ij}&=-L_{ij}&&\text{for $i>j$.}
\end{align*}
\end{theorem}
This theorem, in fact, implies all the obvious constraints between
Seifert form and monodromy, such as the fact that $S=L-L^t$ is
preserved by the local monodromies, as well as the relation
$LH_*(h(\infty))=L^t$, discussed below in the proof of this theorem.

\section{The topological model}\label{sec:alldim}
In this section we will prove the results of sections
\ref{sec:homological} and \ref{sec:monodromy}.

Assume that $f^{-1}(c)$ has just isolated singularities. We first
describe a topological model for the set $N(c)=f^{-1}(D^2_\epsilon(c))$
of Definition \ref{def:vanishing}.

Let $\hat F(c)$ be the result of removing the interiors of the Milnor
fibers $F_i(c)$, $i=1,\dots,t_c$, from $F$.  Topologically, $\hat
F(c)$ results by gluing half-open collars on the regular boundary
components of $F^{ns}(c)$. The local monodromy map $h(c)\colon F\to
F$ can therefore be taken to be the identity on $\hat F(c)$.

\begin{figure}[htbp]
\centerline{\epsfxsize.3\hsize\epsffile{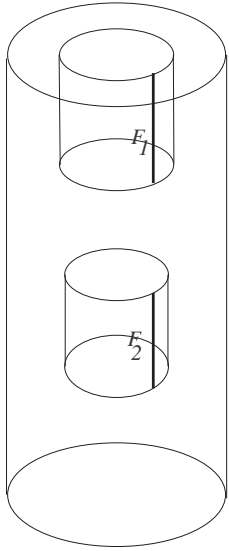}}
\caption{Schematic picture of $N_0$.  Here $F_1$ is a Milnor 
  fiber at infinity, $F_2$ a finite Milnor fiber.}\label{fig:1}
\end{figure}

Define 
\begin{equation*}
  N_0:=(\hat F(c)\times D^2)\cup_\psi \pmt F_{h(c)}\times I,
\end{equation*}
where $\mt F_{h(c)}$ is the mapping torus for the
local monodromy map $h(c)\colon F\to F$ and $\psi$ is the embedding
\begin{equation*}
  \psi\colon \hat F(c)\times S^1=\mt \hat F(c)_{h(c)} \to\pmt
F_{h(c)}\times\{0\}\subset \pmt F_{h(c)}\times I,
\end{equation*}
(See Fig.\ \ref{fig:1}).
The boundary of $N_0$ is the disjoint union of $(t_c+1)$ pieces:
\begin{align*}
   \partial N_0=\bigl(\pmt F_{h(c)}\times \{1\}\bigr)\,\cup\,
 &\bigcup_{i=1}^{t_c}\,\bigl(\pmt
 {F_i(c)}_{h(c)}\cup(\partial F_i(c)\times D^2)\bigr)\cr
 = 
 \qquad\partial_0 N_0\quad\cup\,&\bigcup_{i=1}^{t_c}
 \quad\partial_iN_0
 \qquad\text{(notation)}.
\end{align*}
When $F_i(c)$ is a finite Milnor fiber, that is $s_{c+1}\le i \le t_c$, 
\begin{equation*}
  \partial_iN_0=\pmt {F_i(c)}_{h(c)}\cup
(\partial F_i(c)\times D^2)\cong S^{2n-1}
\end{equation*}
is a standard picture of the sphere with its Milnor fibration for the
link of the singularity in question (see e.g., \cite{neumann-rudolph}).  Let
$N$ be the result of pasting balls $D^{2n}$ to $N_0$ along these
spheres $\partial_iN_0$ for $i=s_{c+1},\dots,t_c$.
\begin{proposition} \label{prop:model} $N$ is a partial compactification of
  $N(c)$ in the following sense: $N(c)$ is homeomorphic\footnote{We
    will prove homeomorphism, but using standard angle straightening
    arguments, cf \cite{hirsch}, one can get a diffeomorphism.}  to
  the result of removing from $N$ all boundary components except the
  ``outer'' boundary component $\partial_0N=\mt F_{h(c)}$.
  \end{proposition}
\begin{proof}
  If $D_i$ be a small enough ball around the $i$-th singularity of
  $f^{-1}(c)$ for $i=s_{c+1},\dots, t_c$ then $f^{-1}(c)$ is transverse to
  each $\partial D_i$ and also to $\partial D^{2n}_r(0)$. By
  compactness, there exists $\epsilon$ so that $f^{-1}(c')$ is also
  transverse to each of these spheres for $|c'-c|\le\epsilon$. Let
  $D_0=D^{2n}_r(0)-\bigcup_{i=1}^{s_c} \interior D_i$ and
  $X=f^{-1}(D^2_\epsilon(c))\cap D_0$.  Then $f|X$ is a submersion of
  a compact manifold with boundary, so by Ehresmann's theorem (see,
  e.g., \cite{hirsch}) it is a locally trivial fibration. Since it is a
  fibration over a disk it is a trivial fibration, so $X\cong
  F^{ns}(c)\times D^2_\epsilon(c)$.  This gives the embedding
  $\phi\colon F^{ns}(c)\times D^2_\epsilon(c)\hookrightarrow \C^n$
  used in the definition of Milnor fibers at infinity.  We can extend
  $\phi$ to the collars outside the regular boundary components of
  $F^{ns}(c)$ (by definition of regular boundary components) to get
  $\phi\colon\hat F(c)\times D^2_\epsilon(c)\hookrightarrow \C^n$
  compatible with the map $f$ and with the trivial structure of $f$
  outside regular boundary components of $F^{ns}(c)$.
  
Let $A=D^2_\epsilon(c)-\interior D^2_{\epsilon/2}(c)$.  Then it is
clear that $X\cup f^{-1}A$ is diffeomorphic to $N_0$ so we will
identify $N_0$ with this subset $X\cup f^{-1}A$ of $N(c)$. The closure
of $N(c)-N_0$ consists of $t_c$ pieces, of which the last $t_c-s_c$
are ``Milnor disk'' neighborhoods of the singularities of $f^{-1}(c)$.
Gluing these back in to $N_0$ gives an embedding $N\to N(c)$, the
closure of whose complement consists of pieces $Y_i$ attached at the
boundary components $\partial_iN:=\pmt {F_i(c)}_{h(c)}\cup(\partial
F_i(c)\times D^2)$ of $N$ for $i=1,\dots,s_c$.  If we show each $Y_i$
is homeomorphic to a half-open collar on $\partial_iN$, then adding
$Y_i$ to $N$ has the same effect up to homeomorphism as removing
$\partial_iN$, so the proof is complete.
  
  To see $Y_i$ is a collar we can use a standard vector-field
  argument. Since $f^{-1}(c)$ is transverse to large spheres about
  $0$, we can find a vector-field $w$ in a neighborhood of any point
  of $f^{-1}(c)$ outside $D^{2n}_r(0)$ so that $w$ is tangent to
  fibers of $f$ and has radially outward component of magnitude $1$.
  Gluing these local $w$'s by a partition of unity, we can find a
  vector-field $w$ which is defined on all $Y_i$, is zero off a
  neighborhood of $f^{-1}(c)\cap Y_i$, is tangent to fibers of $f$, and
  has radially outward component of magnitude $1$ on $f^{-1}(c)\cap
  Y_i$ and of magnitude at most $1$ elsewhere. We can also assume $w$
  is non-zero on the part $\partial F_i(c)\times D^2_\epsilon$ of
  $\partial Y_i$.
  
  Let $v_0$ be the inward radial vector-field $v_0(x,y)=-(x,y)$ on
  $D^2_\epsilon$. Again, by gluing local choices by a partition of
  unity, we can find a vector-field $v$ on $Y_i$ whose image under $f$
  is $v_0$, which is tangent on the part $\partial F_i(c)\times
  D^2_\epsilon$ of $\partial Y_i$, and which has globally bounded
  magnitude.
  
  The sum $v+w$ is then a vector-field on $Y_i$ whose flow-lines all
  lead in backwards time to $\partial Y_i$ and intersect $\partial
  Y_i$ transversally, and whose forward flow lines continue for
  infinite time.  Integrating the vector-field from $\partial Y_i$ thus
  gives a homeomorphism of $Y_i$ with $\partial Y_i\times[0,\infty)$,
  completing the proof.
\end{proof}

\begin{lemma}
  If $N_0$ is constructed as for the above proposition then
  $V_q(c)=\ker\bigl(H_q(F)\to H_q(N_0)\bigr)$ for $q=1,\dots,2n-3$ and
  $V_q(c)=0$ otherwise.  Moreover, $H_q(F)\to H_q(N_0)$ is surjective
  for $q\ne 2n-1$.
\end{lemma}
\begin{proof}
  In the previous proof we identified $N_0$ with a subset of $N(c)$ in
  such a way that $N(c)$ differs from $N_0$ by closing some $S^{2n-1}$
  boundary components by disks and adding collars to some other
  boundary components. It follows that the homology of $N(c)$ and
  $N_0$ differ only in degree $2n-1$. Since $V_q(c)$ vanishes if $q$
  is not in the range $1,\dots,2n-3$, the lemma follows from Theorem
  \ref{th:van}.
\end{proof}
\begin{proof}[Proof of Theorem \ref{th:homology}]
The above lemma implies that $V_q(c)\cong
H_{q+1}(N_0,F)$ by an isomorphism that fits in a commutative diagram:
\begin{equation*}
\begin{CD}
0 @>>> H_{q+1}(N_0,F) @>>> H_q(F) @>>> H_q(N_0) @>>> 0\\
@. @VV\cong V @| @| \\
0 @>>> V_q(c) @>>> H_q(F) @>>> H_q(N_0) @>>> 0
\end{CD}
\end{equation*}

We identify $N_0$ with a subset of $N(c)$ as in the proof of
Proposition \ref{prop:model}. Thus $f$ maps $N_0$ to the disk
$D^2_\epsilon(c)$.  Moreover, $N_0$ is the union of the outer shell
$N_{out}:=f^{-1}(A)$, where $A=D^2_\epsilon(c)-\interior
D^2_{\epsilon/2}(c)$, and an inner core
$N_{inn}=\overline{N_0-N_{out}}$ isomorphic to $\hat F(c)\times
D^2_{\epsilon/2}(c)$.

Express the disk $D^2_{\epsilon}(c)$ as the union of two half-disks
$D^2_-$ and $D^2_+$ by cutting along a diameter.  Let $N_-$ and $N_+$
be the parts of  $N_{out}$ that lie over $D_-$ and
$D_+$ and put $N_1:=N_-\cup N_{inn}$.
We have 
\begin{equation*}
  N_0=N_+\,\cup\,N_1,
\end{equation*}
with
\begin{align*}
    N_{+}&\cong F\times I\times I ,\\
  N_1&\cong F\times I\times I \cup_{\hat F(c)\times I\times\{0\}} \hat
  F(c)\times D^2.
\end{align*}
Thus $N_1$ has $F\times I\times I$ as a deformation retract and the
pair $(N_+,N_+\cap N_1)$ has its intersection with $F\times I\times
\{0\}$ (isomorphic to $\bigl(F\times I, (F\times\partial I)\cup \hat
F(c)\times I\bigr)$) as a deformation retract. In particular, we see
that each of the following maps induces an isomorphism in homology,
since they are, respectively, a homotopy equivalence, an excision map,
and a homotopy equivalence:
\begin{align*}
  (N_0,F)&\hookrightarrow (N_0,N_1)\\
  (N_+,N_+\cap N_1)&\hookrightarrow (N_0, N_1)\\
  \bigl(F\times I, (F\times\partial I)\cup (\hat F(c)\times
  I)\bigr)&\hookrightarrow (N_+, N_+\cap N_1).
\end{align*}
Since
\begin{equation*}
\bigl(F\times I, (F\times\partial I)\cup (\hat F(c)\times
  I)\bigr)=
    (F,\hat F(c))\times (I,\partial I),
\end{equation*}
we get a homology isomorphism
\begin{equation*}
  H_{q+1}(N_0,F) \cong H_{q+1}\bigl((F,\hat F(c))\times
  (I,\partial I)\bigr).
\end{equation*}
The K\"unneth theorem thus gives
\begin{equation*}
H_{q+1}(N_0,F)\cong H_q(F,\hat F(c))\otimes H_1(I,\partial I)
= H_q(F,\hat F(c)).
\end{equation*}
Summarizing, we have an isomorphism 
\begin{equation*}
  H_q(F,\hat F(c))\cong V_q(c).
\end{equation*}

The composition $ H_q(F,\hat F(c))\stackrel{\cong}{\longrightarrow}
V_q(c)\to H_q(F)$ is the variation map (up to sign). Indeed, for the
isomorphism $H_q(F,\hat F(c))\to H_{q+1}(N_0,N_1)$, a relative cycle
$C$ is taken to a relative cycle $C\times I$, mapped by $C\times I\to
F\times I\times\{0\}\subset N_+\subset N_0$ (we are identifying $N_+$
with $F\times I\times I$). We are interested in the boundary of this
cycle as a cycle for $H_q(N_1)\cong H_q(F)$.  When we retract $N_1$ to
$F$, the subset of $N_+\cap N_1$ given by $\partial(F\times
I\times\{0\})$ maps to $F$ by the identity on one component and by $h$
on the other.  The resulting cycle in $F$ thus represents $\pm\var([C])$.

Theorem \ref{th:homology} now follows because $H_q(F,\hat F(c))\cong
\bigoplus_{i=1}^{t_c}H_q(F_i,\partial F_i)$ by excision.
\end{proof}
\begin{proof}[Proof of Theorem \ref{th:monodromy}]
  The first part of Theorem \ref{th:monodromy} just says that the
  image $\Im(H_q(h(c))\colon H_q(F)\to H_q(F))$ is contained in
  $V_q(c)$, which is part of Theorem 1.4 of \cite{neumann-norbury} (it
  is also proved in section 2 of \cite{dimca-nemethi}).  The second
  part of Theorem \ref{th:monodromy} is immediate from the above proof
  of Theorem \ref{th:homology}.
\end{proof}

\begin{proof}[Proof of Theorem \ref{th:seifert}]
  We first recall from \cite{neumann-rudolph} (see also
  \cite{neumann-inv}) how the Seifert linking form for the regular
  link at infinity can be defined on $H_q(F^{co})$.  Let $D^2$ be a
  large disk in $\C$ which contains all $c\in \C$ for which
  $f^{-1}(c)$ is either singular or fails to be M-tame at infinity (in
  the sense of \cite{nemethi-zaharia}, see also
  \cite{tibar-Wallpaper}; there are finitely many such $c$ and they
  include all irregular values of $f$).  Then there is a radius $R$ so
  that for any $r\ge R$ the boundary of the disk $D^{2n}_r(0)$
  intersects all fibers $f^{-1}(t)$ with $t\in \partial D^2$
  transversally.  Then $D:=f^{-1}(D^2)\cap D^{2n}_r(0)$ is
  homeomorphic to $D^{2n}$. The embedding of $F^{co}\in \partial D$ as
  $f^{-1}(t)\cap D$ with $t\in \partial D^2$ gives a Seifert surface
  for the regular link at infinity.  Let $F^{co}_+$ be a neighboring
  copy of $F^{co}$, obtained by replacing $t$ by a nearby point $t_+$
  of $\partial D^2$.  If $a$ is a cycle for homology $H_q(F^{co})$,
  let $a_+$ be a copy of the cycle in $F^{co}_+$.  The Seifert form is
  the form
  \begin{equation*}
    H_q(F^{co})\otimes H_{2n-2-q}(F^{co})\to \Z,\quad
    a\otimes b\mapsto \ell(a_+,b),
  \end{equation*}
  where $\ell$ is linking number in $S^{2n-1}=\partial D$.  It can be
  computed (up to a sign which depends on conventions; following
  \cite{durfee-var} and \cite{neumann-rudolph} the sign is $(-1)^{q+1}$)
  by letting $a_+$ and $b$ bound chains $A_+$ and $B$ in $D$ and
  taking the intersection number $A_+\cdot B$.
  
  We now choose our base point $c_0$ for which $f^{-1}(c_0)$ is our
  ``standard'' regular fiber to be the above point $t$, so we have
  paths $\gamma_c$ (as chosen before Theorem \ref{th:homology}) from
  $t$ to the irregular values $c$.  We can assume these paths run in
  the disk $D^2$.
  
  Suppose that we have a homology class $[a]$ in the image of the map
  $v_i(c)\colon H_*(F_i(c),\partial F_i(c))\to H_*(F)$.  Here, we will
  consider, for the moment, $F$ to be a regular fiber $f^{-1}(c')$
  with $c'\in\partial D^2_\epsilon(c)$, so $F$ is on $\partial N(c)$.
  We write the cycle $a$ as $\var(\alpha)$ with $\alpha$ a relative
  cycle in $(F_i,\partial F_i)$.  By transporting $F_i$ around the
  circle $\partial D^2_\epsilon(c)$ we obtain a map of $\alpha\times
  I$ to $\mt{F_i}_h\subset \partial N(c)$.  The boundary of this chain
  $\alpha\times I$ consists of the union of $\alpha$, $h(\alpha)$, and
  $\partial \alpha\times I$.  The part $\partial \alpha\times I$
  bounds a mapping of $\partial \alpha\times D^2_\epsilon(c)$ mapping
  to $\partial F_i\times D^2_\epsilon(c)\subset N(c)$, so we can glue
  this to $\alpha\times I$ to get a chain $A_0$ in $N(c)$ with
  boundary $\partial A_0$ representing $var(\alpha)=a$ mapping to $F$.
  
  If we want $a$ be a cycle for homology of our standard fiber
  $F=f^{-1}(t)$ we glue onto the above $A_0$ a copy of $a\times I$
  mapping into $f^{-1}(\gamma_c)$.  We call the resulting chain $A$.
  Note that $A_0$ lies completely in the ``shell''
  $\pmt {F_i}_h\cup(\partial F_i\times D^2_\epsilon(c)\subset N_0\subset
  N(c)$.  We can also construct $A_0$ in a smaller shell obtained by
  replacing $\epsilon$ by $\epsilon/2$ and removing a thin collar from
  $\partial F_i$.  We denote the version of $A$ constructed this way
  by $A^{thin}$.
  
  Suppose now the two homology classes $[a],[b]$ are in the image of
  the map $v_i(c)\colon H_*(F_i(c),\partial F_i(c))\to H_*(F)$, where
  $F$ is now our standard regular fiber.  We can assume they both lie
  in $F^{co}$, since $F$ retracts to $F^{co}$. We can then make $b$
  bound a cycle $B$ as above.  We can also make $a_+$ bound a cycle
  $A_+^{thin}$ constructed as above but using a path $(\gamma_c)_+$
  running parallel to $\gamma_c$ to a point on $\partial
  D^2_{\epsilon/2}(c)$.  This path runs through a point $c'_+$ next
  to $c$ on $\partial D^2_\epsilon(c)$. The chains $A_+^{thin}$ and
  $B$ intersect only in the fiber $f^{-1}(c'_+)$ and the intersection
  number $A_+^{thin}\cdot B$ is, up to sign, the intersection number
  in $f^{-1}(c')$ of $a$ and $\beta$. (See Fig.~2.)
\begin{figure}[htbp]
\centerline{\hfill \epsfxsize.3\hsize\epsffile{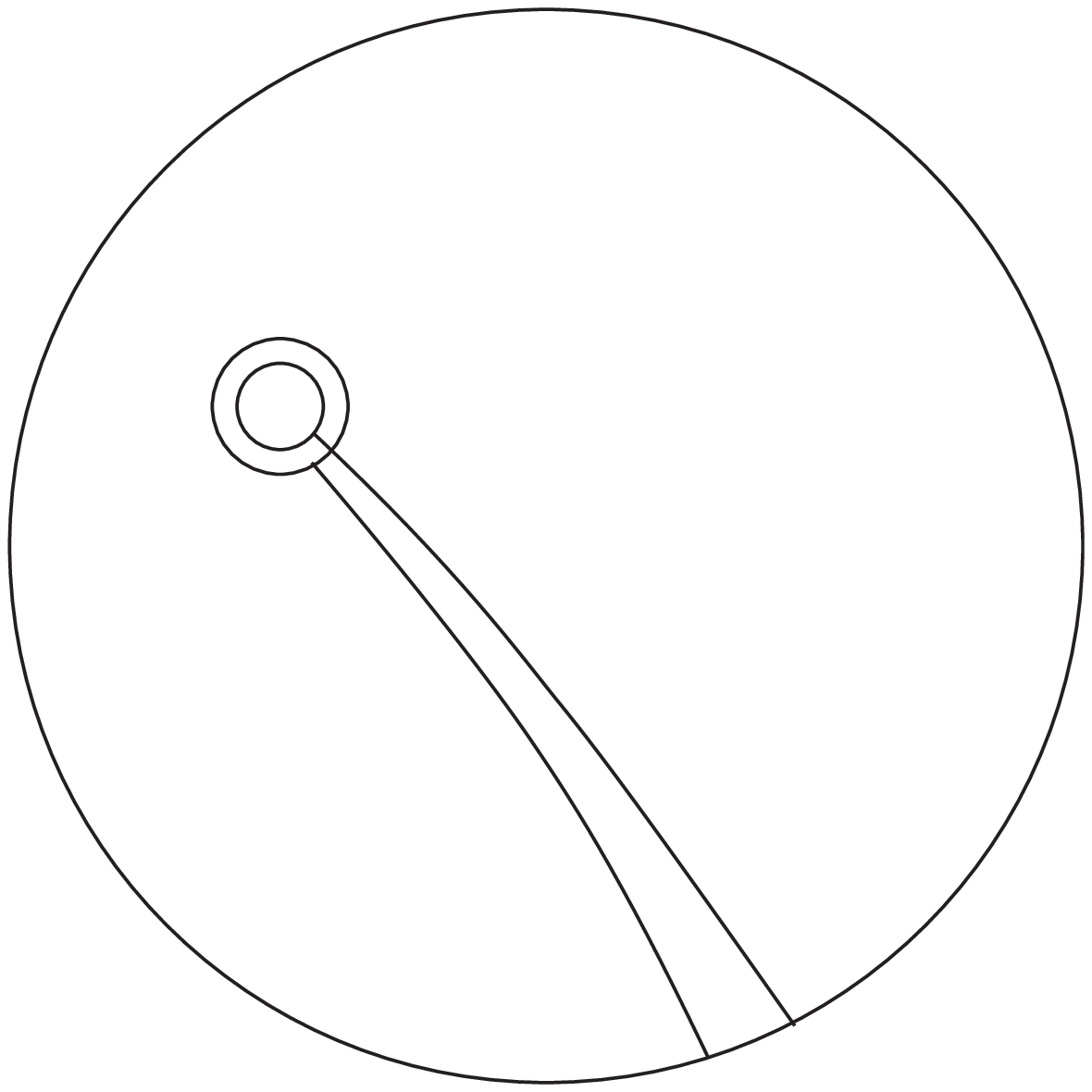}\hfill 
\epsfxsize.3\hsize\epsffile{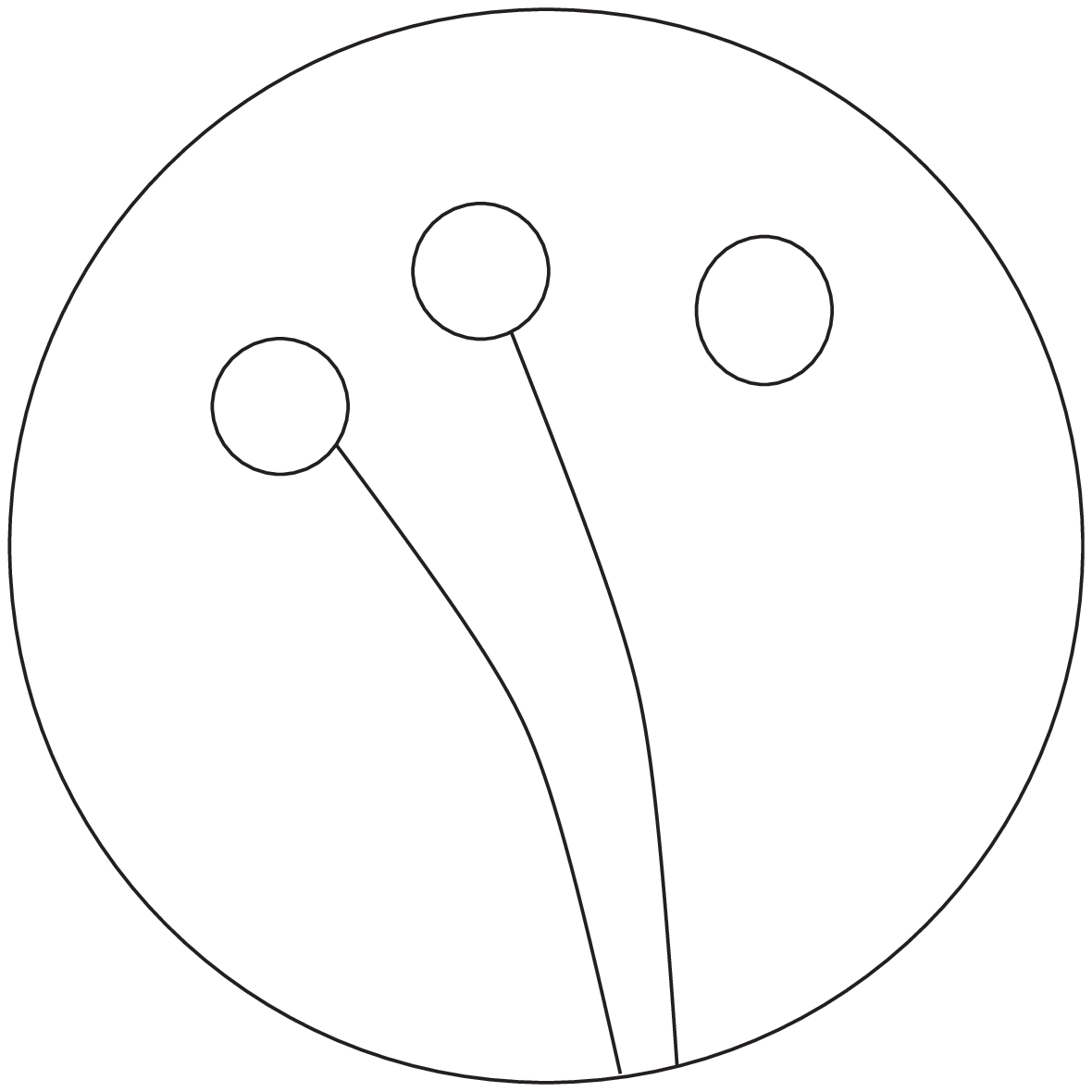} \hfill}
\centerline{\hfill \textsc{Figure 2}\qquad\hfill\qquad \textsc{Figure 3}\hfill}
\end{figure}\stepcounter{figure}\stepcounter{figure}
With standard sign conventions, the sign is in fact $+1$ (this is most
easily checked by using the standard formulae $L^t=LH$ and $S=L-L^t$
relating Seifert form $L$, intersection form $S$ and monodromy $H$,
since for a knot the relationship to be proved is $L(I-H)=S$). This
proves the claim of Theorem \ref{th:seifert} about the diagonal blocks
of the Seifert form.  The claim about vanishing of appropriate
off-diagonal blocks is the same as the corresponding proof in
\cite{neumann-rudolph}, as suggested by Fig.~3.\end{proof}
\begin{proof}[Proof of Theorem \ref{th:relations}]
  Let $M$ be a 0-codimensional submanifold (with boundary) of
  the sphere $S^m$ and suppose $M$ is fibered over the circle $S^1$
  with fiber $F$.  Then we can define a Seifert form $L$ and
  homological monodromy $H$ on the homology of $F$ as for fibered
  links and the obvious geometric relation
  $L(x,Hy)=\ell(x_+,Hy)=\ell(x,y_+)=L^t(y,x)$ can be written in matrix
  form as $LH=L^t$ (this is a well-known relation in the case of
  fibered links, see, e.g., \cite{durfee-var}).  If we apply this in
  the situation of Theorem \ref{th:relations} it gives the equation
  $$LH_*(h(c_k))H_*(h(c_{k-1}))\dots H_*(h(c_1))=L^t.$$
  An inductive
  argument, which we omit, shows that this equation is equivalent to
  the collection of equations of the theorem. 
\end{proof}

\section{Polynomials in dimension 2}
In the remainder of this paper we describe results specific to
dimension $2$.  We assume that $f\colon\C^2\to \C$ has only isolated
singularities.

The regular link at infinity determines and is determined by a certain
fibered multilink \cite{neumann-inv}, which we will call the
\emph{fundamental multilink}. By \cite{neumann-irreg} the link at
infinity of any irregular fiber is obtained by splicing additional
links onto this multilink; we can call these the \emph{splice
  components at infinity}.  They may decompose further as splices of
fibered and non-fibered parts, and the fibered parts are called the
\emph{fibered splice components at infinity}. 
\begin{theorem}\label{th:dim2}
  For $n=2$ the Milnor fibers at infinity and their monodromy maps
  arise as the fibers and monodromy of the fibered splice components at
  infinity.
\end{theorem}
We describe below how the fibered splice components at infinity are
determined by the splice diagram for the link. This only depends on
the topology of the link (see \cite{eisenbud-neumann}), so we have the
important corollary:
\begin{corollary}
  For $n=2$ the Milnor fibers at infinity and their monodromy maps are
  completely determined by the link at infinity of the irregular
  fibers they belong to (and are effectively computable from their
  splice diagrams, as described in \cite{eisenbud-neumann}).\qed
\end{corollary}

The \emph{splice diagram} of the link at infinity of a complex affine
plane curve (\cite{neumann-inv}), which from the point of view of
classical algebraic geometry is simply an encoding of the Puiseux tree
at infinity, also encodes the splice decomposition of the link at
infinity (\cite{eisenbud-neumann}).  It is a weighted tree with some
leaves drawn as arrowheads to stand for link components of the link.
Splicing links corresponds to gluing such diagrams at arrowheads.
Conversely, disconnecting a splice diagram by cutting an edge and
drawing two arrowheads on the resulting ends corresponds to the
inverse operation of splice decomposition.  Any link obtainable via
repeated splice decomposition is called a ``splice component.'' Splice
components are thus represented by connected subgraphs of the splice
diagram.

The \emph{fibered splice components at infinity} are the splice
components corresponding to the maximal connected subgraphs of the
splice diagram having only negative vertex linking weights. The
example below will clarify this.

\vspace{6pt} Before proving Theorem \ref{th:dim2}, we illustrate it
using the example of the ``Brian\c con polynomial''
$$f(x,y)=x^2(1+xy)^4+3x(1+xy)^3+(3-\frac83x)(1+xy)^2 -4(1+xy)+y.$$
This polynomial was shown to have no finite singularities by Brian\c
con, see \cite{ACL} where it is also shown that all fibers of $f$ are
connected. It has two irregular fibers (over $0$ and $-16/9$). The
Jordan normal forms for action in homology of the monodromy generators
$h(0)$ and $h(-16/9)$ were computed by Artal-Bartolo, Cassou-Nogues,
and Dimca \cite{artal-cassou-dimca}. Dimca and Nemethi
\cite{dimca-nemethi} computed these with respect to a common basis of
homology, thus determining the complex monodromy representation for
this example.  We will show how the splice diagrams make these
computations ``routine'' and give the geometric monodromy rather than
just the action on homology.  However, as the example will make clear,
our approach still falls short of achieving our goal of a practical
complete algorithmic description of the topology.

The splice diagrams for the links at infinity of the fibers of $f$
were computed in \cite{ACL}. The regular splice diagram is as
follows, where we have included the linking weights (also called
multiplicity weights) at vertices in parentheses: 
$$
\objectmargin{0pt}\spreaddiagramrows{-5pt}
\spreaddiagramcolumns{12pt}
\diagram
\\
&\overtag\Circ{(0)}{16pt}
\ulto\dlto\ddline^(.25){2}\rline^(.25){-3}^(.75){1}&\Dot
\rline^(.25){1}^(.75){-1}&\overtag\Circ{(2)}{16pt}
\ddline^(.25){2}\rline^(.25){1}^(.75){-7}
&\overtag\Circ{(3)}{16pt}\ddline^(.25){3}\rto&\\ \\
&\lefttag\Circ{(0)}{6pt}&&\lefttag\Circ{(1)}{6pt}&\lefttag\Circ{(1)}{6pt}
\enddiagram
.
$$
The fact that vertices with zero linking weights occur is equivalent to
each of the following two facts (see \cite{eisenbud-neumann},
\cite{neumann-inv} and \cite{neumann-irreg}):
\begin{itemize}
\item the regular link at infinity is not a fibered link;
\item $f$ has fibers that are irregular at infinity.
\end{itemize}
We write the link as the splice of the part with zero linking
weights and a fibered multilink:
$$
\objectmargin{0pt}\spreaddiagramrows{-5pt}
\spreaddiagramcolumns{12pt}
\diagram
\\
&\overtag\Circ{(0)}{16pt}
\ulto\dlto\ddline^(.25){2}\rto^(.25){-3}&\righttag~{(6)}{3pt}&
\lefttag~{(4)}{3pt}&\Dot\lto_(.25){1}
\rline^(.25){1}^(.75){-1}&\overtag\Circ{(2)}{16pt}
\ddline^(.25){2}\rline^(.25){1}^(.75){-7}
&\overtag\Circ{(3)}{16pt}\ddline^(.25){3}\rto&\\ \\
&\lefttag\Circ{(0)}{6pt}&&&&\lefttag\Circ{(1)}{6pt}&\lefttag\Circ{(1)}{6pt}
\enddiagram
$$
The fibered multilink is the multilink associated with the regular
link at infinity for $f$ as described in \cite{neumann-inv}. We call
it the \emph{fundamental multilink} for $f$. As described there, it is
a fibered multilink which determines and is determined by the regular
link at infinity.  Its fibers are, up to isotopy, the regular fibers
of $f$ over the points of a large circle in $\C$ so its fibration
gives the monodromy at infinity for $f$ (which is the product of the
local monodromies around the irregular fibers). We will return to this
later and first examine the irregular fibers.

The Brian\c con polynomial has irregular fibers $f^{-1}(0)$ and
$f^{-1}(-16/9)$.  The link at infinity for $f^{-1}(0)$ has splice diagram
$$
\objectmargin{0pt}\spreaddiagramrows{-5pt}
\spreaddiagramcolumns{12pt} \diagram &\overtag\Circ{(0)}{16pt}
\lto\ddline^(.25){2}^(.75){-2}\rline^(.25){-3}^(.75){1}&\Dot
\rline^(.25){1}^(.75){-1}&\overtag\Circ{(2)}{16pt}
\ddline^(.25){2}\rline^(.25){1}^(.75){-7}
&\overtag\Circ{(3)}{16pt}\ddline^(.25){3}\rto&\\ \\
&\righttag\Circ{(-1)}{6pt}\ulto\dlto&&\lefttag\Circ{(1)}{6pt}&
\lefttag\Circ{(1)}{6pt}\\
& \enddiagram
$$
It follows that $f^{-1}(0)$ has Euler characteristic $-2$. Since it has
$4$ boundary components, it is a four-punctured sphere. 

\begin{figure}[htbp]
\centerline{\epsfxsize.25\hsize\epsffile{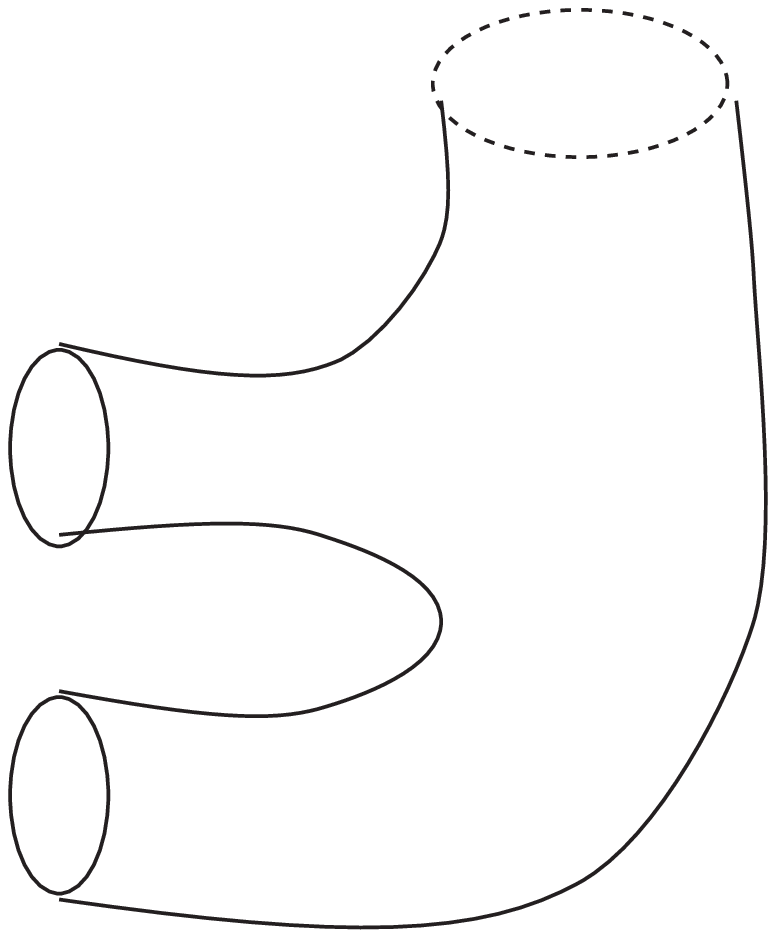}}
\caption{}\label{fig:5}
\end{figure}

We express the link at infinity of this irregular fiber as the
splice of the parts with positive, zero, and negative linking weights
respectively:
$$
\objectmargin{0pt}\spreaddiagramrows{-5pt}
\spreaddiagramcolumns{12pt} \diagram &\overtag\Circ{(0)}{16pt}
\lto\dto^(.4){2}\rto^(.25){-3}&\righttag~{(6)}{2pt}&\lefttag~{(4)}{2pt}
&\lto_(.25){1}\Dot
\rline^(.25){1}^(.75){-1}&\overtag\Circ{(2)}{16pt}
\ddline^(.25){2}\rline^(.25){1}^(.75){-7}
&\overtag\Circ{(3)}{16pt}\ddline^(.25){3}\rto&\\ 
&\undertag~{(2)}{2pt}\\
&\overtag~{(3)}{8pt}&&&&\lefttag\Circ{(1)}{6pt}&
\lefttag\Circ{(1)}{6pt}\\
&\righttag\Circ{(-1)}{6pt}\ulto\dlto\uto_(.4){-2} \\
&\enddiagram
$$
The part with positive linking weights is always the fundamental
multilink (see \cite{neumann-irreg}). We call the parts with negative
total linking weights (in this case there is just one) the
\emph{fibered splice components at infinity}.  So the fibered splice
component at infinity is given by the splice diagram:
$$
\objectmargin{0pt}\spreaddiagramrows{-5pt}
\spreaddiagramcolumns{12pt} \diagram\\&\overtag~{(3)}{8pt}\\ \\
&\righttag\Circ{(-1)}{6pt}\ulto\dlto\uuto_(.25){-2} \\
&\enddiagram
$$
Its fiber has Euler characteristic $-1$ and $3$ boundary
components, so it is a thrice-punctured disk (see Fig.~\ref{fig:5}).
If we remove the boundary component marked with a dashed curve
(corresponding to the arrowhead where the fibered splice component at
infinity splices to the rest of the irregular spliced diagram) we
obtain the Milnor fiber at infinity for this irregular fiber. This
would be $F_1(0)$ in our earlier notation, but there is
just the one Milnor fiber at infinity for $f^{-1}(0)$, so we 
call it $F(0)$.  

The local monodromy for this Milnor fiber is given by the monodromy of
the fibered splice component at infinity.  The book
\cite{eisenbud-neumann} describes how to compute this monodromy from
the splice diagram.  We need the monodromy in which the boundary of
$F(0)$ (consisting of the two circles at the left of the figure) is
fixed. By Theorem 13.5 of \cite{eisenbud-neumann} it is the result of
doing a Dehn twist on an annulus parallel to each of these boundary
components.  It follows that the variation map
\begin{equation*}
  \var\colon \Z\cong H_1(F(0),\partial F(0))\to H_1(F(0))
\end{equation*}
takes a generator of $H_1(F(0),\partial F(0))$ to the homology class
of the difference of the two boundary components.  Following the
orientation conventions of \cite{eisenbud-neumann} the bilinear form
occurring in Theorem \ref{th:seifert} (given by equation
(\ref{eq:seifert})) thus has matrix $(2)$.

\begin{figure}[htbp]
\centerline{\epsfxsize.5\hsize\epsffile{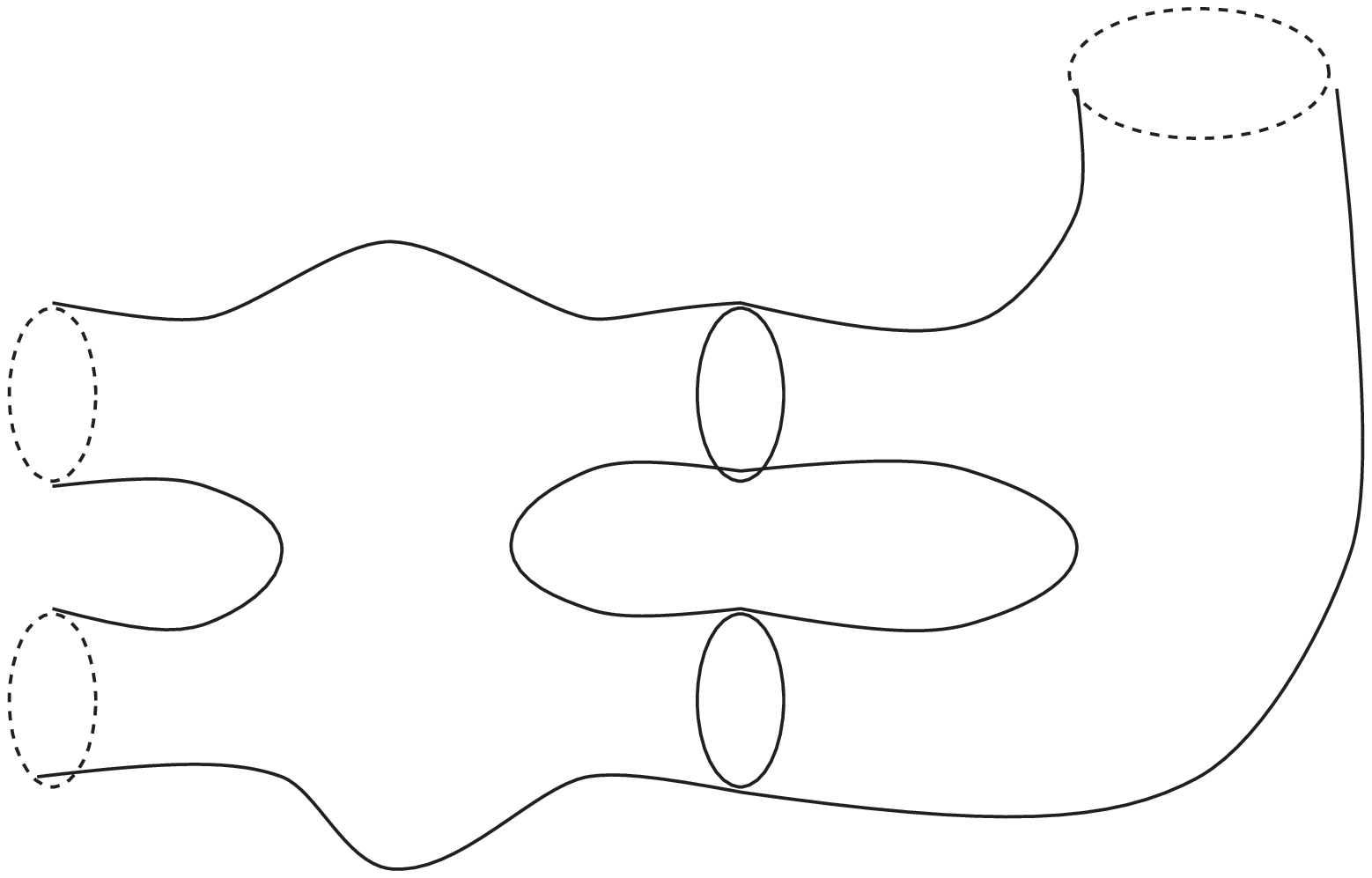}}
\caption{}\label{fig:6}
\end{figure}

The regular fiber $F$ is the result of gluing the Milnor fiber at infinity
to the irregular fiber (Fig.~\ref{fig:6}).
By Theorem \ref{th:homology}, the subgroup $V_1(0)\subset H_1(F)$ of
vanishing cycles for the fiber $f^{-1}(0)$ is generated by the
difference of the two separating curves in this figure.  This is fixed 
by the local monodromy, so the monodromy matrix is
$(1)$.

We will number our irregular values $c_1=-16/9, c_2=0$, since this is
the ordering Dimca and Nemethi use in \cite{dimca-nemethi}.  So, in
the notation of equations (\ref{eq:H}) and (\ref{eq:L}) in section
\ref{sec:monodromy} we have
\def\aaa{p}\def\bbb{q}\def\ccc{r}
\begin{gather}
  \label{eq:c2}
  L_2=(2),\quad h_2=(1), \\
\label{eq:hc2}
H_1(h(c_2))=
  \begin{pmatrix}
    I&0\\h_{21}&h_2
  \end{pmatrix}=
\left(
  \begin{array}{rrr|r}
1&0&0&0\\
0&1&0&0\\
0&0&1&0\\
\hline
\aaa&\bbb&\ccc&1
  \end{array}\right)
\end{gather}
with $h_{21}=(\aaa,\bbb,\ccc)$ still to be determined.

\begin{figure}[htbp]
\centerline{\epsfxsize.5\hsize\epsffile{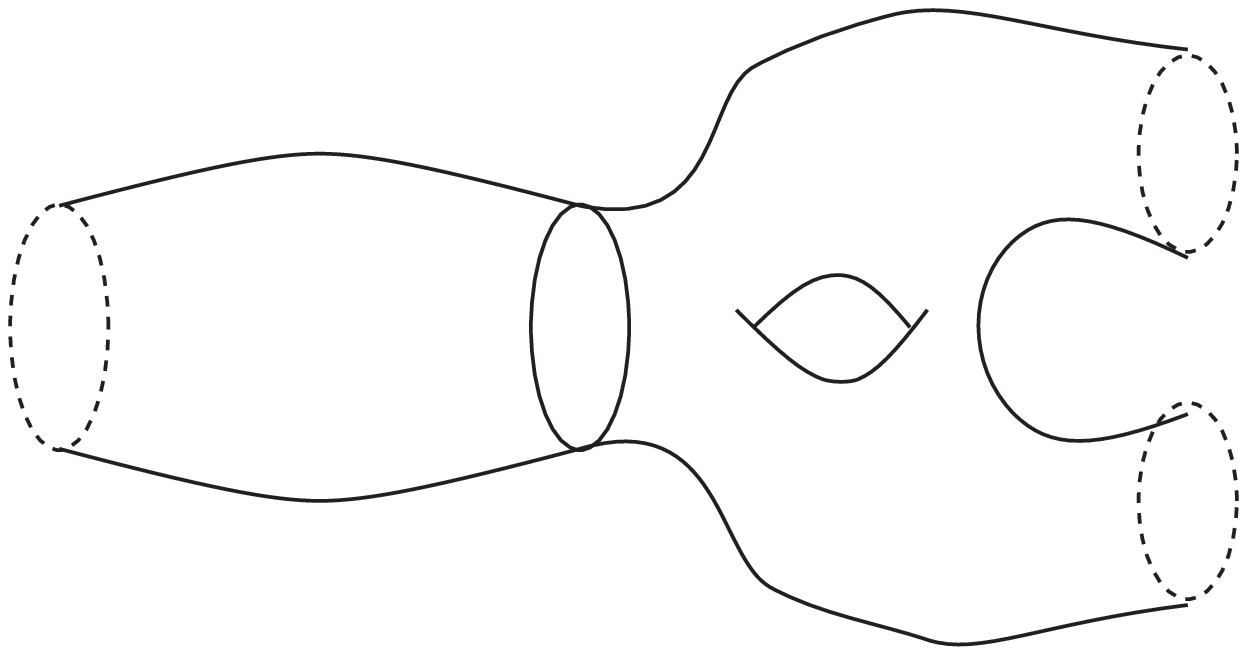}}
\caption{}\label{fig:8}
\end{figure}

We now do a similar analysis for the irregular fiber
$f^{-1}(-16/9)$. The splice diagram for the link at infinity of this
fiber is:
$$
\objectmargin{0pt}\spreaddiagramrows{-5pt}
\spreaddiagramcolumns{12pt}
\diagram\\
&\lto\overtag\Circ{(-6)}{16pt}\rline^(.25){-15}^(.75){1}\ddline^(.25){2}
&\overtag\Circ{(0)}{16pt}\ddline^(.25){2}\rline^(.25){-3}^(.75){1}&\Dot
\rline^(.25){1}^(.75){-1}&\overtag\Circ{(2)}{16pt}
\ddline^(.25){2}\rline^(.25){1}^(.75){-7}
&\overtag\Circ{(3)}{16pt}\ddline^(.25){3}\rto&\\ \\
&\lefttag\Circ{(-3)}{6pt}&\lefttag\Circ{(0)}{6pt}&&\lefttag\Circ{(1)}{6pt}
&\lefttag\Circ{(1)}{6pt}
\enddiagram
$$
so the irregular fiber $f^{-1}(-16/9)$ has Euler characteristic $0$
and is thus an annulus.  Moreover, the fibered splice component at
infinity for this irregular link at infinity has diagram:
$$
\objectmargin{0pt}\spreaddiagramrows{-5pt}
\spreaddiagramcolumns{12pt}
\diagram
&\lto\overtag\Circ{(-6)}{16pt}\rto^(.25){-15}\ddline^(.25){2}
&\righttag~{(12)}{3pt}
\\ \\
&\lefttag\Circ{(-3)}{6pt}
\enddiagram
$$
The fiber of this fibered multilink is a thrice-punctured torus
(right hand piece of Fig.~\ref{fig:8}) which glues to the irregular fiber as in
Fig.~\ref{fig:8} 
to give a copy of the regular fiber $F$. We call this Milnor fiber at
infinity $F(-16/9)$. The local monodromy on it is isotopic to an order
$6$ map (because of the linking weight $-6$) and it exchanges the two
boundary components at the right (the circles corresponding to a
single edge of the splice diagram of a fibered multilink are always
permuted transitively by the monodromy for the fibration).  The local
monodromy on $F$ is thus also isotopic to this order $6$ map, since
$F$ and $F(-16/9)$ just differ by a collar.

Again we can use this description to compute the local monodromy and
the block $L_1$ of the Seifert form.  We describe this in greater
detail later, but a quick approach is to note that there are exactly
two different order $6$ transformations of the surface in question
with given action on the boundary, differing only in orientation. The
correct orientation can be deduced from the boundary twist
computations in \cite{eisenbud-neumann} or by means of the equivariant
signature computation of Theorem 5.3 of \cite{neumann-splice} (as
generalized in section 6 of that paper).  In any case, with respect to 
a suitable basis of homology, the answer is:
\begin{gather}\label{eq:c1}
  L_1=
  \begin{pmatrix}
    1&1&0\\0&1&0\\0&0&0
  \end{pmatrix},\quad
h_1=\left(
\begin{array}{rrr}
  0&-1&0\\1&1&0\\0&-1&-1
\end{array}\right)\\ \label{eq:hc1}
H_1(h(c_1))=
  \begin{pmatrix}
    h_1&h_{12}\\0&I
  \end{pmatrix}=
\left(
  \begin{array}{rrr|r}
0&-1&0&a\\
1&1&0&b\\
0&-1&-1&c\\
\hline
0&0&0&1
  \end{array}\right)
\end{gather}
with $h_{12}=(a,b,c)^t$ still to be determined.  

At this point we can write down the Seifert matrix $L$ as
\begin{equation}
  \label{eq:seif}
L=
  \begin{pmatrix}
    L_1&0\\L_{21}&L_2
  \end{pmatrix}=
\left(
  \begin{array}{rrr|r}
1&1&0&0\\
0&1&0&0\\
0&0&0&0\\
\hline
x&y&z&2
  \end{array}\right)
\end{equation}
with $L_{21}=(x,y,z)$ still to be determined.  However, applying the
relations of Theorem \ref{th:relations} gives:
\begin{equation*}
  a+b=x,~b=y,~0=z,\quad 2\aaa=-x,~2\bbb=-y,~2\ccc=-z,  
\end{equation*}
so in fact we just have three unknown integers in equations
(\ref{eq:hc2}), (\ref{eq:hc1}), (\ref{eq:seif}), namely
$\aaa,$ $\bbb,$ $c$, and the others are then determined as
\begin{equation}
  \label{eq:consolidate}
  \ccc=z=0,~x=-2\aaa,~y=-2\bbb,~a=2\bbb-2\aaa,~b=-2\bbb.
\end{equation}

\begin{figure}[htbp]
\centerline{\epsfxsize.7\hsize\epsffile{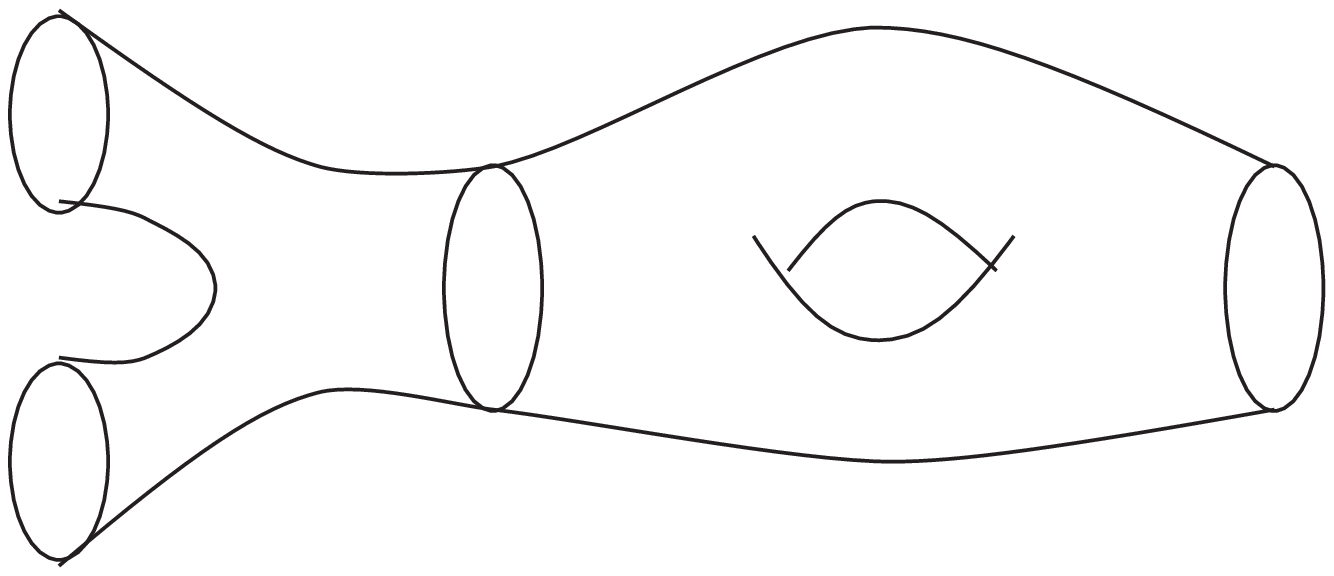}}
\caption{}\label{fig:9}
\end{figure}

The product of the two local monodromy maps is the monodromy at
infinity for the regular fiber.  This monodromy at infinity is 
the monodromy of the fundamental multilink (by the definition of
this multilink in \cite{neumann-inv}), so it can be computed
from the splice diagram for this multilink.  In our
particular case that splice diagram (as an unrooted diagram) is
$$
\objectmargin{0pt}\spreaddiagramrows{-5pt}
\spreaddiagramcolumns{12pt}
\diagram\\
\lefttag~{(4)}{3pt}
&\overtag\Circ{(2)}{16pt}\lto_(.25){-1}
\ddline^(.25){2}\rline^(.25){1}^(.75){-7}
&\overtag\Circ{(3)}{16pt}\ddline^(.25){3}\rto&\\ \\
&\lefttag\Circ{(1)}{6pt}&\lefttag\Circ{(1)}{6pt}
\enddiagram
$$
The fiber of this fibered multilink decomposes according to the
splice components determined by the two nodes of this diagram as in
Fig.~\ref{fig:9},
and the monodromy restricted to the left part is isotopic to a map of
order $2$ and restricted to the right part is isotopic to a map of
order $3$.  The sixth power of this monodromy gives a single Dehn
twist on the joining circle by Theorem 13.1 of
\cite{eisenbud-neumann}. 

This monodromy map in $H_1(F)$ has eigenvalues $-1,1,e^{\pm2\pi i/3}$,
so its characteristic polynomial is $(t+1)(t^3-1)=t^4+t^3-t-1$.  On
the other hand, equations (\ref{eq:hc2}), (\ref{eq:hc1}), and
(\ref{eq:consolidate}) show the monodromy map is
\begin{gather*}
H_1(h(c_2))H_1(h(c_1))=
  \begin{pmatrix}
  1&0&0&0\\
0&1&0&0\\
0&0&1&0\\
\aaa&\bbb&0&1
\end{pmatrix}
\begin{pmatrix}
  0&-1&0&-2\aaa+2\bbb\\
1&1&0&-2\bbb\\
0&-1&-1&c\\
0&0&0&1
\end{pmatrix}
\\ \qquad\qquad=
\begin{pmatrix}
 0&-1&0&-2\aaa+2\bbb\\ 
1&1&0&-2\bbb\\
0&-1&-1&c\\
\bbb&-\aaa+\bbb&0&-2\aaa^2+2\aaa\bbb-2\bbb^2+1
\end{pmatrix}.
\end{gather*}
This has characteristic polynomial
$$(t-1)(t^3+2(\aaa^2-\aaa\bbb+\bbb^2-1)(t^2-t)-1),$$ so
$$\aaa^2-\aaa\bbb+\bbb^2-1=0.$$
This equation has six solutions: 
$$(\aaa,\bbb)=\pm(1,1),~\pm(1,0),~\pm(0,1)$$
with corresponding values
$$(a,b)=\pm(0,-2),~\pm(-2,0),~\pm(2,-2).$$
But, by changing our choice
of basis on $H_1(F(-16/9))$ by powers of the local monodromy we cycle
through these six possibilities, so they are all equivalent.  Choosing
$(\aaa,\bbb)=(1,0)$ gives the conclusion:
\begin{align*}
H_1(h(c_2))&=
\begin{pmatrix}
  1&0&0&0\\
0&1&0&0\\
0&0&1&0\\
1&0&0&1
\end{pmatrix}\\
H_1(h(c_1))&=
\begin{pmatrix}
  0&-1&0&-2\\
1&1&0&0\\
0&-1&-1&c\\
0&0&0&1
\end{pmatrix}
\\ 
L&=
\begin{pmatrix}
  1&1&0&0\\
0&1&0&0\\
0&0&0&0\\
-2&0&0&2
\end{pmatrix}
\end{align*}

Here $c$ is still undetermined, and
we know no way of finding it with our current methods. However, one
calculates easily that the isomorphism type of the monodromy
representation over $\C$ depends only on the vanishing or not of
$3c+2$, which is non-vanishing since it is not divisible by $3$. An
equivalent non-vanishing issue arose in the computation of complex
monodromy in \cite{dimca-nemethi} and was resolved by a more
complicated argument.

We have described what can be read directly and easily from the splice
diagrams.  To complete the information about the global monodromy
takes more work, since we must identify the three different pictures
of the regular fiber $F$ of Figs.~\ref{fig:6}, \ref{fig:8}, and
\ref{fig:9} to fully understand the global picture.  It is easy to see
that the boundary component at the right of Fig.~\ref{fig:9}
corresponds to the one at the left in Fig.~\ref{fig:8} and one of the
ones at the left in Fig.~\ref{fig:6}.  The issue is to determine how
the two circles in Fig.~\ref{fig:6} lie with respect to the order $6$
map of Fig.~\ref{fig:8}. This would also determine how they lie with
respect to the separating circle of Fig.~\ref{fig:9}. The homology
information gives strong hints, but no obvious complete answer. In the
next section we give a conjectural answer.

\section{A tentative picture of the Brian\c con topology}

Let $F$ denote the three-punctured torus. We will
describe explicit maps $h_1$ and $h_2$ of $F$ that satisfy all the
properties of the local monodromy maps $h(c_1)$ and $h(c_2)$ for the
Brian\c con polynomial that were computed in the previous section.
Namely, $h_1$, $h_2$, and $h_2h_1$ are conjugate in the group of
orientation preserving diffeomorphisms of $F$ to the maps $h(c_1)$,
$h(c_2)$, and $h(\infty)$ of the previous section.

\begin{figure}[htbp]
  \begin{center}
    \epsfxsize.4\hsize\epsffile{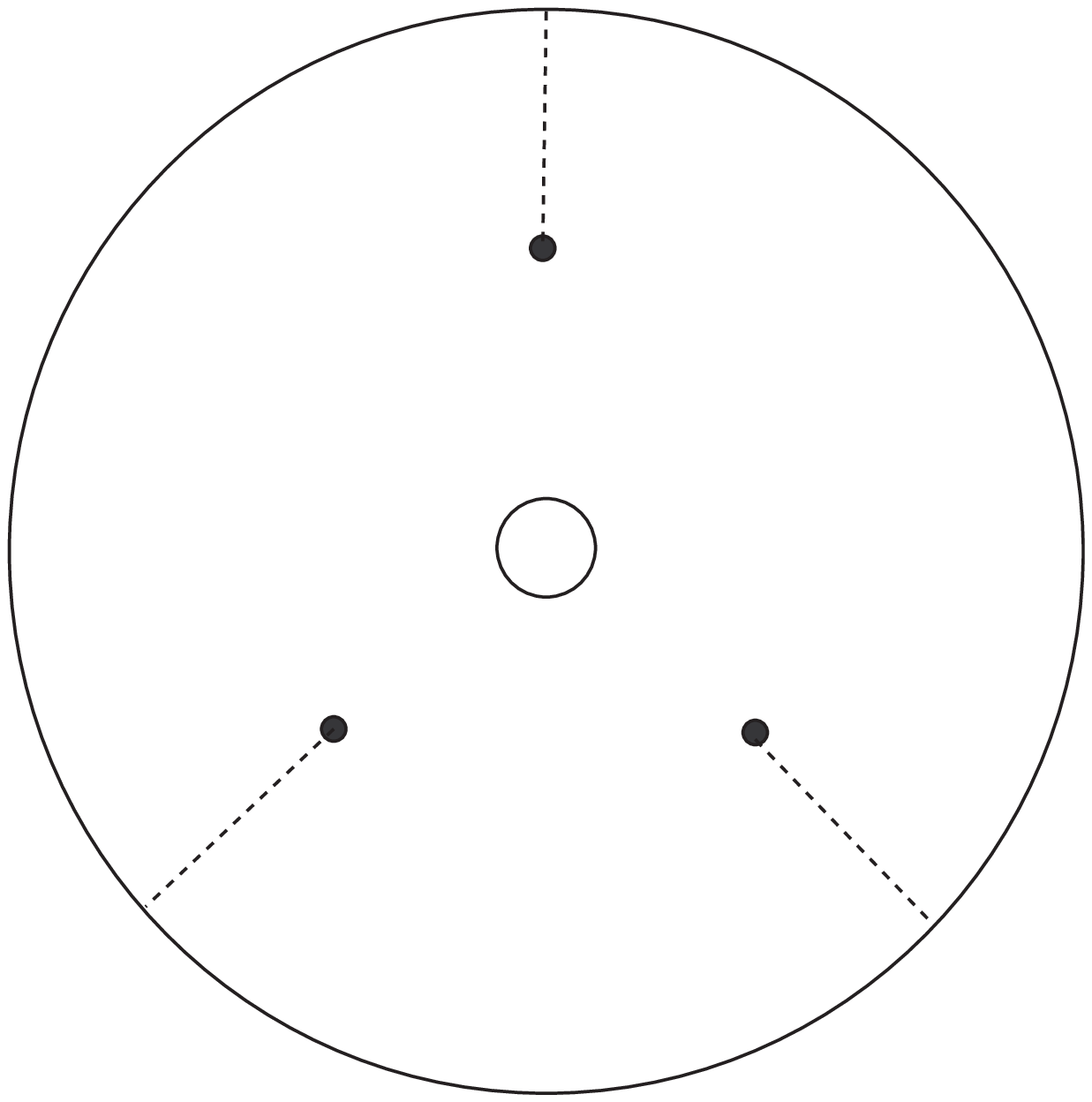}
  \end{center}
    \caption{{$F$ as a $2$-fold cover; the dashed lines represent branch cuts}}
    \label{fig:10}
\end{figure}

We will represent $F$ as the 2-fold cover of a punctured disk branched
at three points, described by branch cuts as in Fig.~\ref{fig:10}. 
The inner boundary component is thus covered by two boundary
components of $F$ and the outer boundary component is double covered
by one boundary component of $F$.

Our map $h_1$ will be the order $6$ map which rotates the picture by
one-third of a turn clockwise and exchanges the two branches.

\begin{figure}[htbp]
  \begin{center}
    \epsfxsize.4\hsize\epsffile{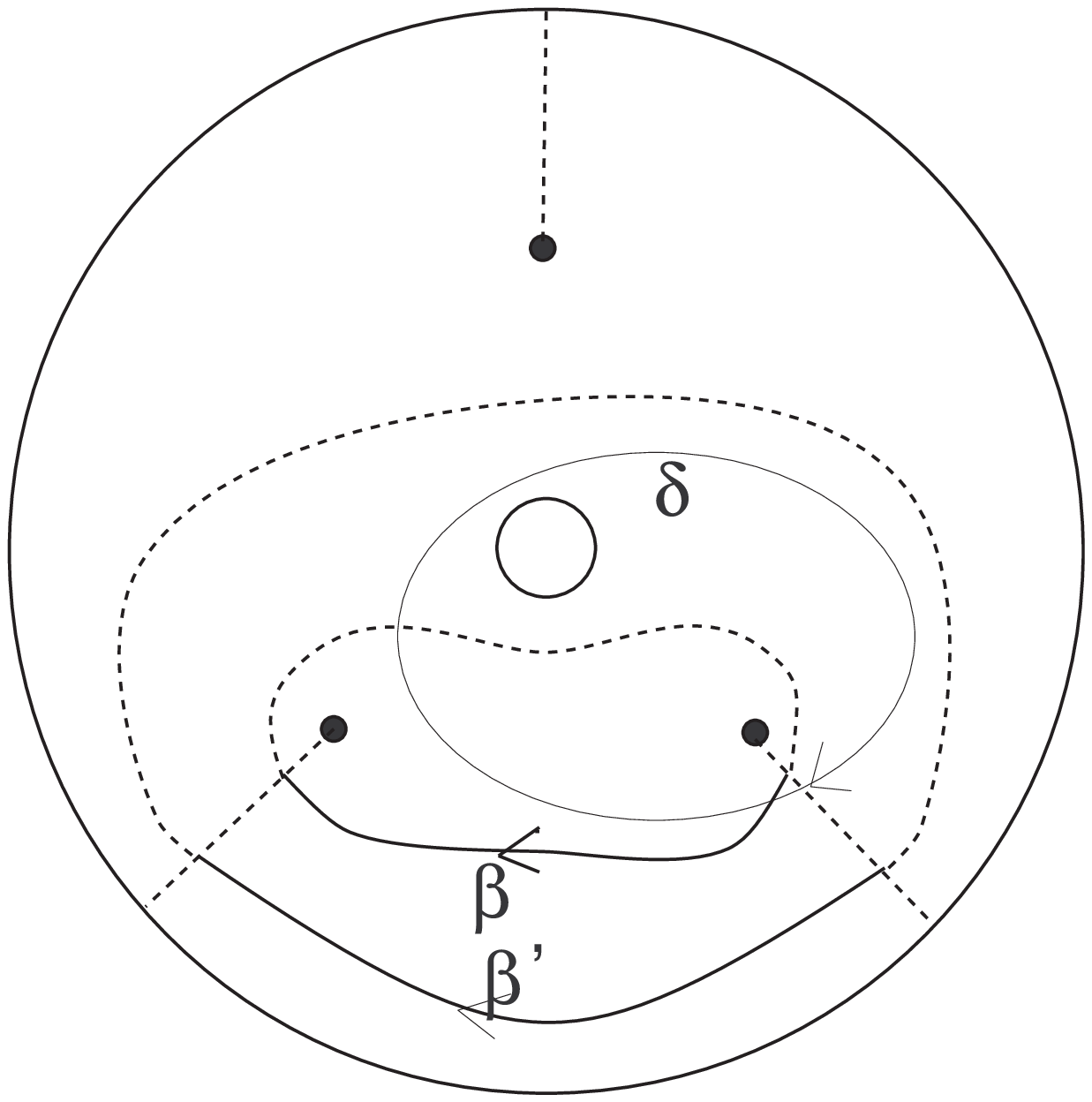}
  \end{center}
    \caption{The dashed portions of curves lie on the ``back'' branch
      of $F$; $\delta$ is double covered in $F$}
    \label{fig:13}
\end{figure}

Consider now the curves on $F$ labeled $\beta$, $\beta'$, $\delta$,
as in Fig.~\ref{fig:13}.
The curves $\beta$ and $\beta'$ cut $F$ into two pieces of genus zero
as in Fig.~\ref{fig:6}. We therefore define $h_2$ to be the product of
the Dehn twists on these two curves.  The curve $\delta$ cuts $F$ into
two pieces as in Fig.~\ref{fig:9}. By drawing a careful picture one
finds that the effect of $h_2$ on $h_1(\delta)$ is to take it to
$\delta$, so $h_2h_1$ takes $\delta$ to itself.  In fact,
$h_2h_1$ is conjugate to the map $h(\infty)$, with $\delta$ playing
the role of the separating curve in Fig.~\ref{fig:9}.  This can be
seen by drawing careful pictures, but it is also forced by the fact
that $\delta$ is mapped to itself together with the homology
computation below.

\begin{figure}[htbp]
  \begin{center}
    \epsfxsize.4\hsize\epsffile{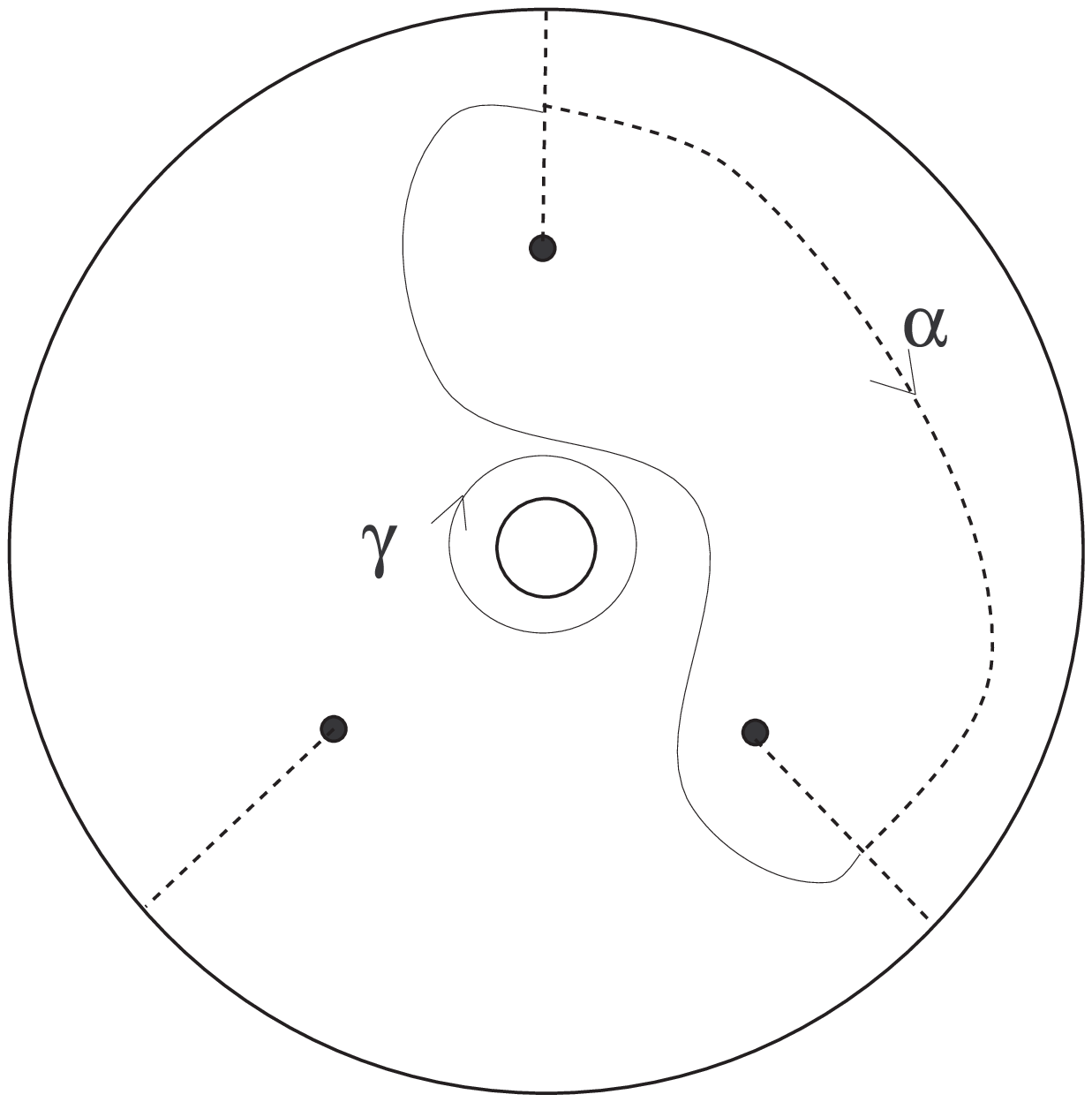}
  \end{center}
    \caption{}
    \label{fig:11}
\end{figure}

Now let $\alpha$ and $\gamma$ be the curves of Fig.~\ref{fig:11} and
let $\gamma'$ be the curve represented by a dashed circle in place of
$\gamma$ (this is the image of $\gamma$ under $h_1$).  The
intersection number $\alpha\cdot\beta$ is $1$, so $\alpha$ and $\beta$
form a basis for the homology $H_1(\overline F)$ of the closed torus
obtained by filling the punctures of $F$. A basis for the homology of
$F$ is therefore given by $\alpha$, $\beta$, $\gamma$, $\gamma'$.
However, computing the images of the variation maps for $h_1$ and
$h_2$ leads to the basis 
$$\alpha,~ \beta,~ \gamma-\gamma';\quad \beta+\beta'=2\beta+\gamma',$$
which we therefore use instead.  With respect to this basis it is
easily checked that the actions of $h_2$ and $h_1$ on homology are by
the matrices
$$
\begin{pmatrix}
  1&0&0&0\\
0&1&0&0\\
0&0&1&0\\
1&0&0&1
\end{pmatrix}\quad\text{and}\quad
\begin{pmatrix}
  0&-1&0&-2\\
1&1&0&0\\
0&-1&-1&-1\\
0&0&0&1
\end{pmatrix}.$$
These are the matrices 
for $H_1(h(c_2))$ and $H_1(h(c_1))$ of the previous
section with $c=-1$.

As already stated, we do not know if the above is really the
appropriate monodromy.  The computation of the previous section
implies that the two curves of Fig.~\ref{fig:6} are
$\beta+r(\gamma-\gamma')$ and $\beta+r(\gamma-\gamma')\pm\gamma'$ in
homology for some $r\in\Z$ (and it then follows easily that the $c$ in
$H_1(h(c_2))$ is odd: namely $c=-4r-2\pm1$).  There are many pairs of
disjoint simple closed curves that satisfy this, but the fact that
$h(c_2)h(c_1)$ has to fix a separating closed curve and be isotopic to
order two and three maps on the resulting pieces of $F$ seems a very
strong constraint, and may well eliminate most or all other
possibilities.

\section{Proof of Theorem \ref{th:dim2}}

We first need to recall some basics about multilinks and splicing.
See \cite{eisenbud-neumann} for details.

\vspace{3pt}
For the moment, by a ``link'' we will understand a pair $(\Sigma, L)$
consisting of an oriented submanifold $L$ of dimension 1 in a
3-dimensional homology sphere $\Sigma$. It is a ``knot'' if $L$
consists of a single closed curve. Our homology sphere $\Sigma$ will
always be $S^3$ in applications in this paper, but the discussion of
splicing is easier without this restriction.  

The \emph{link exterior} for a link $(\Sigma,L)$ is the manifold with
boundary $\Sigma-\interior N(L)$, where $N(L)$ is a (small) closed
regular neighborhood of $L$ in $\Sigma$.

A \emph{multilink} is a link $(\Sigma, L)$ with an integer
``multiplicity'' $m(K)$ assigned to each component $K$ of $L$, with
the convention that reversing the orientation of a component $K$ and
simultaneously changing the sign of $m(K)$ gives the same multilink.
In other words, the multilink structure is a given by a 1-cycle $m$
supported on $L$.  Equivalently, and more conveniently, it is given by
the cohomology class $\mu\in H^1(\Sigma-L;\Z)$ whose value on a
$1$-cycle $c$ is the linking number $\ell(m,c)$.  A Seifert surface
for the multilink is a map of a compact oriented surface $S$ to
$\Sigma$ which maps $\partial S$ to $L$, is an embedding on
$S-\partial S$, and, considered as a $2$-chain, has the above
$1$-cycle $m$ as boundary.  If $N(L)$ is a regular neighborhood of $L$
that intersects the Seifert surface $S$ in a collar on $\partial S$
then the surface $S-S\cap\interior N(L)$ in the link exterior
$\Sigma-\interior N(L)$ is also called a Seifert surface.

If $(\Sigma_1,K_1)$ and
$(\Sigma_2,K_2)$ are knots, we form the \emph{splice} 
\def\splice#1#2#3#4{#1 \frac{}{#2~~~~~~#4}#3}
\begin{equation*}
\Sigma=\splice{\Sigma_1}{K_1}{\Sigma_2}{K_2}
\end{equation*}
by pasting together link exteriors of each knot as follows:
\begin{equation*}
\Sigma=(\Sigma_1-\interior N_1)\cup_\partial(\Sigma_2-\interior N_2),
\end{equation*}
where the pasting along boundaries $\partial N_1$ and $\partial N_2$
is done so as to match a meridian of $K_1$ to a longitude of $K_2$ and
meridian of $K_2$ with longitude of $K_1$.  A simple homology
calculation shows $\Sigma$ is again a homology sphere.

If $K_1$ and $K_2$ are components of links $L_1\subset\Sigma_1$ and
$L_2\subset\Sigma_2$ and $L=(L_1-K_1)\cup(L_2-K_2)$ then we
write
\begin{equation*}
  (\Sigma,L)=\splice{(\Sigma_1,L_1)}{K_1}{(\Sigma_2,L_2)}{K_2},
\end{equation*}
the \emph{splice of $(\Sigma_1,L_1)$ to $(\Sigma_2,L_2)$ along $K_1$
and $K_2$}.  

If $(\Sigma,L)$ has a multilink structure then we get
induced multilink structures on each $(\Sigma_i, L_i)$ by restricting
the cohomology class that defines the multilink structure.  Note that,
even if the multilink $(\Sigma, L)$ is a link (all multiplicities are
$1$), the multiplicities of $K_1$ and $K_2$ will in general be
different from $1$.  Thus decomposing links via splicing leads one
naturally into the realm of multilinks.

If $(\Sigma,L$) results from splicing two links as above then there is
a torus $T^2$ in the link exterior $\Sigma-\interior N(L)$ along which
the splicing occurred.  Conversely, suppose $(\Sigma, L)$ is a link and
$T^2\subset \Sigma-\interior N(L)$ an essential torus (i.e., the
induced mapping $\pi_1(T^2)\to \pi_1(\Sigma-\interior N(L))$ is
injective and $T^2$ is not isotopic to a boundary component of
$\Sigma-\interior N(L)$). Then $(\Sigma,L)$ is the result of a
non-trivial splicing operation along this torus: cutting $\Sigma$
along $T^2$ gives two homology solid tori, and we obtain $\Sigma_1$
and $\Sigma_2$ from $\Sigma$ by replacing each of these homology solid
tori in turn by a genuine solid torus.  

In this situation, if $\Sigma$ is $S^3$ then $\Sigma_1$ and $\Sigma_2$
are also $S^3$.  In the following we will start with a link in $S^3$
and splice decompose it, so we never see homology
spheres other than $S^3$.

\begin{figure}[htbp]
\centerline{\epsfxsize.4\hsize\epsffile{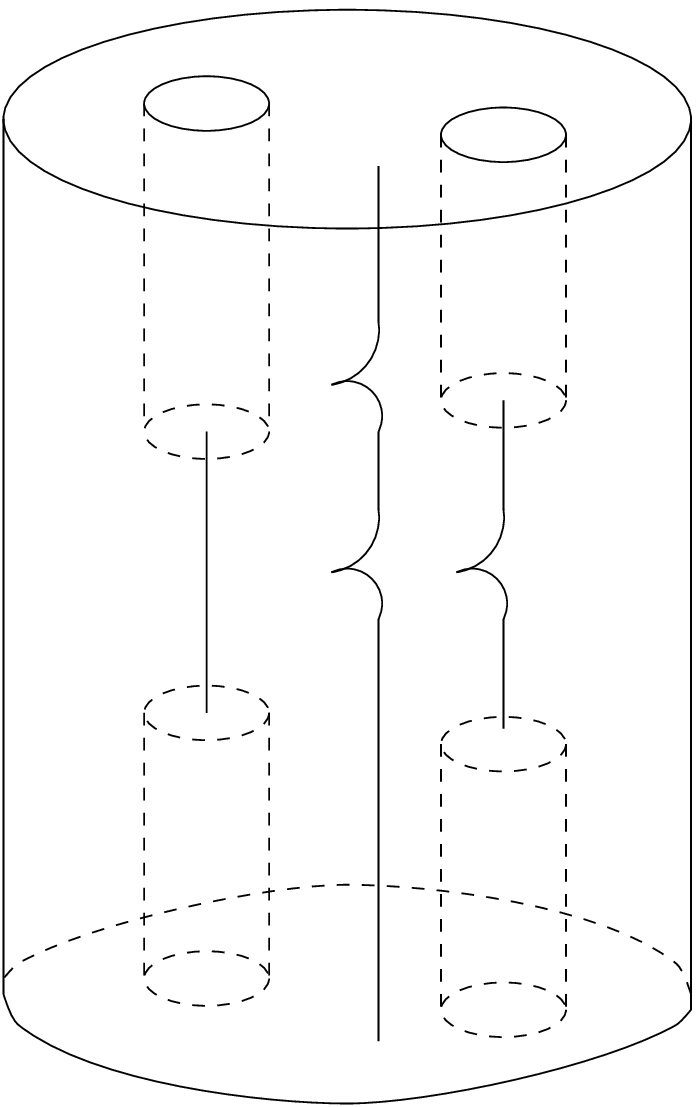}}
\caption{}\label{fig:4}
\end{figure}

\def\D{{\mathbb D}}
We will study the topology of $f\colon\C^2\to \C$ by intersecting
fibers with a large disk in $\C^2$.  In fact our basic topological
object will be a $4$-disk $\D$ obtained as follows: take a large disk
$D^2$ that contains all irregular values of $f$, intersect
$f^{-1}(D^2)$ with a very large $4$-disk in $\C^2$, and then push in
``holes'' around the fibers that are irregular at infinity as in the
models $N$ of Proposition \ref{prop:model}.  See Fig.~\ref{fig:4}.

We will need to do this carefully to confirm the desired properties of
the resulting space.  As in \cite{neumann-inv}, it is convenient to
use a polydisk $D(q,r):=\{(x,y)\in\C^2:|x|\le q,|y|\le r\}$ for our
``very large $4$-disk''.  We recall Lemma 2.1 of \cite{neumann-inv}:

\begin{lemma}
  By a linear change of coordinates we may assume $f(x,y)$ is of
  degree $n$ and of the form
  $f(x,y)=x^n+f_{n-1}(y)x^{n-1}+\dots+f_0(y)$. We choose $s$ so all
  irregular values of $f$ lie in the disk $D^2_s(0)\subset\C$. Then
  for $r$ sufficiently large and $q$ sufficiently large with respect
  to $r$ the fibers $f^{-1}(t)$ for $t\in\partial D^2_s(0)$ intersect
  $\partial D(q,r)$ only in the part $|x|<q,|y|=r$, and do so
  transversely --- in fact they intersect each line $y=y_0$ with
  $|y_0|\ge r$ transversely.
\end{lemma}

We sketch a slight modification of the argument in \cite{neumann-inv}.
The fiber $f^{-1}(t)$ fails to be transverse to the line $y=y_0$ if and
only if $y_0$ is the image of a branch point of the projection
$f^{-1}(t) \to \C$ given by the $y$-coordinate.  If $f^{-1}(t)$ is
reduced (no multiple components) there will be finitely many such
branch points on $f^{-1}(t)$. The locus of such branch points as $t$
varies is an algebraic curve $B$ in $\C^2$ (given by the equation
$\Delta(y,f(x,y))=0$, where $\Delta(y,t)$ is the discriminant of the
polynomial $f(x,y)-t\in\C[y,t][x]$). A fiber $f^{-1}(c)$ is irregular
at infinity if and only if it is not reduced (in which case it has a
component in common with $B$) or if intersection points of $B$ with
nearby fibers $f^{-1}(t)$ move off to infinity as $t$ approaches $c$.
Thus $f^{-1}(\partial D^2_s(0))\cap B$ is compact, and if we choose
$r$ large enough that this compact set lies in the domain $|y|<r$ of
$\C^2$ then $r$ does what is desired.

This proof actually shows more.  Choose $\epsilon$ small enough that
the disks $D^2_\epsilon(c)$ about the irregular values of $f$ are
pairwise disjoint and lie in the interior of $D^2_s(0)$.  Denote
$X=D^2-\bigcup_{c\in\Sigma}\interior D^2_\epsilon(c)$.  Then
\begin{scholium}
With notation as above, if the radius $r$ is sufficiently large then
the fibers $f^{-1}(t)$ for $t\in X$ intersect
each line $y=y_0$ with
  $|y_0|\ge r$ transversely.
\end{scholium}

As described in \cite{neumann-inv}, if $r$ is chosen as in the above
Lemma, then
\begin{equation*}D:=f^{-1}(D^2_s(0))\cap\{|y|\le r\}
\end{equation*}
is a $4$-disk in $\C^2$.  (This can be seen by noting that $\C^2$
results by gluing $f^{-1}(D^2_s(0))\cap\{|y|\ge r\}$ to $D$ along part
of $\partial D$ and then gluing $f^{-1}(\C-\interior D^2_s(0))$ along
the boundary of the result. The first part glued on is a collar
because $f^{-1}(D^2_s(0))\cap\{|y|\ge r\} \to [r,\infty)$,
$(x,y)\mapsto |y|$ is a locally trivial fibration, and the second part
is obviously a collar.)

We now assume $r'$ was chosen sufficiently large that all irregular
fibers $f^{-1}(c), c\in \Sigma$, are transverse to the lines $y=y_0$
with $|y_0|\ge r'$. In particular, the irregular fibers are transverse
to the cylinders $\{|y|=r\}$ with $r\ge r'$. Then assume $\epsilon$
was chosen small enough that all fibers $f^{-1}(t)$ with $t$ in
$\bigcup_{c\in\Sigma}D^2_\epsilon(c)$ are transverse to the cylinders
$\{|y|=r'\}$.  Then $r$ is chosen as in the Scholium above. 

For each $c\in\Sigma$ the set $f^{-1}\bigl(\interior D^2_\epsilon(c)\bigr)
\cap\bigl\{|y|>r'\bigr\}$ consists of components $\left(f^{-1}\bigl(\interior D^2_\epsilon(c)\bigr)
\cap\bigl\{|y|>r'\bigr\}\right)_i$, $i=1,\dots,r_c$ corresponding to
places where $f^{-1}(c)$ is irregular at infinity, and maybe
additional components where $f^{-1}(c)$ is regular at infinity.
Let
\begin{equation*}
  \D=D-\biggl(\bigcup_{c\in\Sigma}\bigcup_{i=1}^{r_c}\left(f^{-1}\bigl(\interior D^2_\epsilon(c)\bigr)
\cap\bigl\{|y|>r'\bigr\}\right)_i\biggr),
\end{equation*}
see Fig.~\ref{fig:4} above.  The argument of Section \ref{sec:alldim}
easily adapts to show that $D$ results topologically by adding a
collar to part of the boundary of $\D$, so $\D$ is homeomorphic to
$D^4$.

We will need names for the parts of the boundary of $\D$.  Denote
\begin{equation*}
E:=f^{-1}\bigl(\partial D^2_s(0)\bigr)\cap \D,\quad S:=\partial
\D-\interior E.
\end{equation*}
Then $S$ is the union of the part 
\begin{equation*}
  S_0:=S\cap\{|y|=r\},
\end{equation*}
and $S-\interior S_0$, which is the union of pieces 
\begin{equation*}
T_i(c)\cong \pmt{F_i(c)}_h\cup\bigl(\partial F_i(c)\times D^2\bigr)   
\end{equation*}
for $c\in\Sigma$ and $i=1,\dots,s_c$ (recall that $\mt{F_i(c)}_h$
denotes the mapping torus of the local monodromy on the Milnor fiber
$F_i(c)$ at infinity). 
\begin{lemma}
  $S$ is a union of solid tori and each $T_i(c)$ is homeomorphic to a
  solid torus.
\end{lemma}
\begin{proof}
  Let $S'=\partial D -\interior E$. That $S'$ is a union of solid tori
  was proved in \cite{neumann-inv}.  The argument (due to L. Rudolph)
  is that for $|y_0|=r$ the intersection $S'\cap\{y=y_0\}$ is
  transverse and the result is a union of disks by the maximum modulus
  principle, since it is equivalent to the set $\{x\in\C:|f(x,y_0)\le
  s\}$. On the other hand, $S'\cong S$ by the argument that identifies
  $D-\interior \D$ with collars on the $T_i(c)$'s, so $S$ is a union of
  solid tori.
  
  The same argument applies to show the sets
  $f^{-1}\bigl(D^2_\epsilon(c)\bigr)\cap\{|y|=r\}$ are unions of solid
  tori for $c\in\Sigma$. But the components of these sets are
  homeomorphic to the $T_i(c)$'s.
\end{proof}

The above decomposition of $\partial \D$ gives splice decompositions of
the links at infinity of the fibers of $f$. The basic fact was
described earlier in this section: if $(S^3,L)$ is a link and we cut
the link exterior $S^3-\interior N(L)$ into pieces along embedded
tori, then this represents $(S^3,L)$ as the result of a splicing
operations.

In particular, the piece $E$ of $\partial \D$ is the exterior of a
splice component $(S^3,L)$, where $L\subset S^3=\partial \D$ is the
link consisting of the cores of the solid tori making up $S$.  This
splice component is the fundamental multilink for $f$ as described in
\cite{neumann-inv}. The fibration of the the exterior $E$ of this
multilink is simply given by the restriction $f|E$.

We shall see that, possibly after minor modification to eliminate
parallel tori, the above splice decomposition of the link at infinity
of an irregular fiber of $f$ is as described in Theorem \ref{th:dim2}.

Note that $f|S_0\colon S_0\to X$ is a fibration of $S_0$ over the
punctured disc $X$ with each fiber a union of circles (isotopic to the
regular link at infinity of $f$).  We can extend this map over the
solid tori $T_i(c)$ to get a Seifert fibration of $S$. But $S$ is a
disjoint union of one or more solid tori and, up to isotopy, the only
Seifert fibrations of a solid torus are the standard
$(p,q)$-fibrations in which the core circle is a fiber and the general
fiber $p$-fold covers this core circle (we do not rule out the
possibility of $(p,q)=(0,1)$, called a ``generalized Seifert
fibration'' in \cite{jankins-neumann}\footnote{One can show it cannot
  occur here, but we do not need this.}).  Thus, each component of
$S-\bigcup_i T_i(c)$ is a solid torus with a collection of thinner
solid tori removed, all or all but one of which run parallel to fibers
of this Seifert fibration, and maybe one running parallel to the core
circle.

The link at infinity of an irregular fiber $f^{-1}(c)$ can be seen as
the intersection of $f^{-1}(c)$ with $\partial \D$. We thus have a splice
decomposition of this link at infinity into:
\begin{itemize}
\item the fundamental multilink,
\item fibered multilinks based on the pieces $T_i(c)$,
\item multilinks with exteriors given by the components of the Seifert
fibered piece $S-\bigcup_i T_i(c)$.
\end{itemize}
The latter will lead to the non-fibered components mentioned in the
theorem.  However, some components of $S-\bigcup_i T_i(c)$ could be of
the form: solid torus minus a thinner solid torus parallel to the core
circle, giving a toral annulus $T^2\times I$.  In this case the two
torus boundary components of this piece are parallel, so to obtain an
irredundant splice decomposition we must omit one of them and absorb
this toral annulus as a collar on an adjacent splice component.

To complete the proof of Theorem \ref{th:dim2} we must show the splice
decomposition we have found is as described in that theorem.  We do
this by examining our construction in terms of a compactification of
$\C^2$.  Since the relationship between the compactification divisor
and the splice diagram is already worked out in detail in
\cite{eisenbud-neumann} and \cite{neumann-irreg}, this then does what
we require.

\def\Div{Y}
We extend the polynomial map $f\colon \C^2\to \C$ to a map $\fbar
\colon Z\to \C P^1$ of a smooth compact complex surface $Z$ to $\C
P^1$. The compactification divisor $\Div:=Z-\C^2$ is a union of smooth
rational curves with dual intersection graph a tree. A component of
$\Div$ on which $\fbar$ is non-constant is called \emph{horizontal}. A
component of $\Div$ on which $\fbar$ is constant is called \emph{finite}
or \emph{infinite} according to the value of $\fbar$ on it.

By blowing up if necessary, we can assume that the only singularities
of fibers of $\fbar$ that occur on $\Div$ are normal crossings between
components of the fiber and components of $\Div$.

We can encode the topology of $\Div$ in the usual way by a plumbing
graph. This is a tree, with vertices corresponding to components of
$\Div$ and edges for intersections between components.  It has a weight at
each vertex to show the self-intersection number of the corresponding
component of $\Div$. We draw arrows at vertices to indicate where fibers
of $\fbar$ intersect $\Div$. The following diagram, which gives a
compactification divisor for the Brian\c con polynomial, uses solid
arrows for the general fiber and dashed respectively dotted arrows for
the two irregular fibers. Note that the three curves on the left could
be blown down if we only wanted a compactification on which $\fbar$ is
well defined; they arose from blowing up to resolve a singularity of
an irregular fiber on $\Div$.
$$\xymatrix@R=12pt@C=12pt@M=0pt@W=0pt@H=0pt{
\overtag\Circ{-2}{9pt}\lineto[rrr]&&&
\overtag\Circ{-1}{9pt}\lineto[rrr]\ar@{-->}[dd]&&&
\overtag\Circ{-3}{9pt}\lineto[rrr]&&&
\overtag\Circ{-1}{9pt}\lineto[rrr]\lineto[dd]
\ar[ddr]\ar[drr]\ar@{.>}[dll]&&&
\overtag\Circ{-2}{9pt}\lineto[rrr]\lineto[dd]&&&
\overtag\Circ{-2}{9pt}\lineto[rrr]\lineto[dd]&&&
\overtag\Circ{-2}{9pt}\lineto[rrr]&&&
\overtag\Circ{-2}{9pt}\\&&&&&&&&&&&&&&&&&&\\
&&&&&&&&&\lefttag\Circ{-2}{6pt}\ar@{.>}[dll]\ar@{.>}[ddl]&&&
\lefttag\Circ{-2}{6pt}&&&
\lefttag\Circ{-2}{6pt}\lineto[dd]\\&&&&&&&&&&&&&&&&&&\\
&&&&&&&&&&&&&&&\lefttag\Circ{-2}{6pt}\lineto[dd]
\\&&&&&&&&&&&&&&&&&&\\
&&&&&&&&&&&&&&&\lefttag\Circ{-1}{6pt}
\ar@{.>}[dll]\ar@{-->}[dd]\ar[drr]
\\&&&&&&&&&&&&&&&&&&\\&&&&&&&&&&&&&&&&&&
}$$

\begin{figure}[htbp]
\centerline{\epsfxsize.9\hsize\epsffile{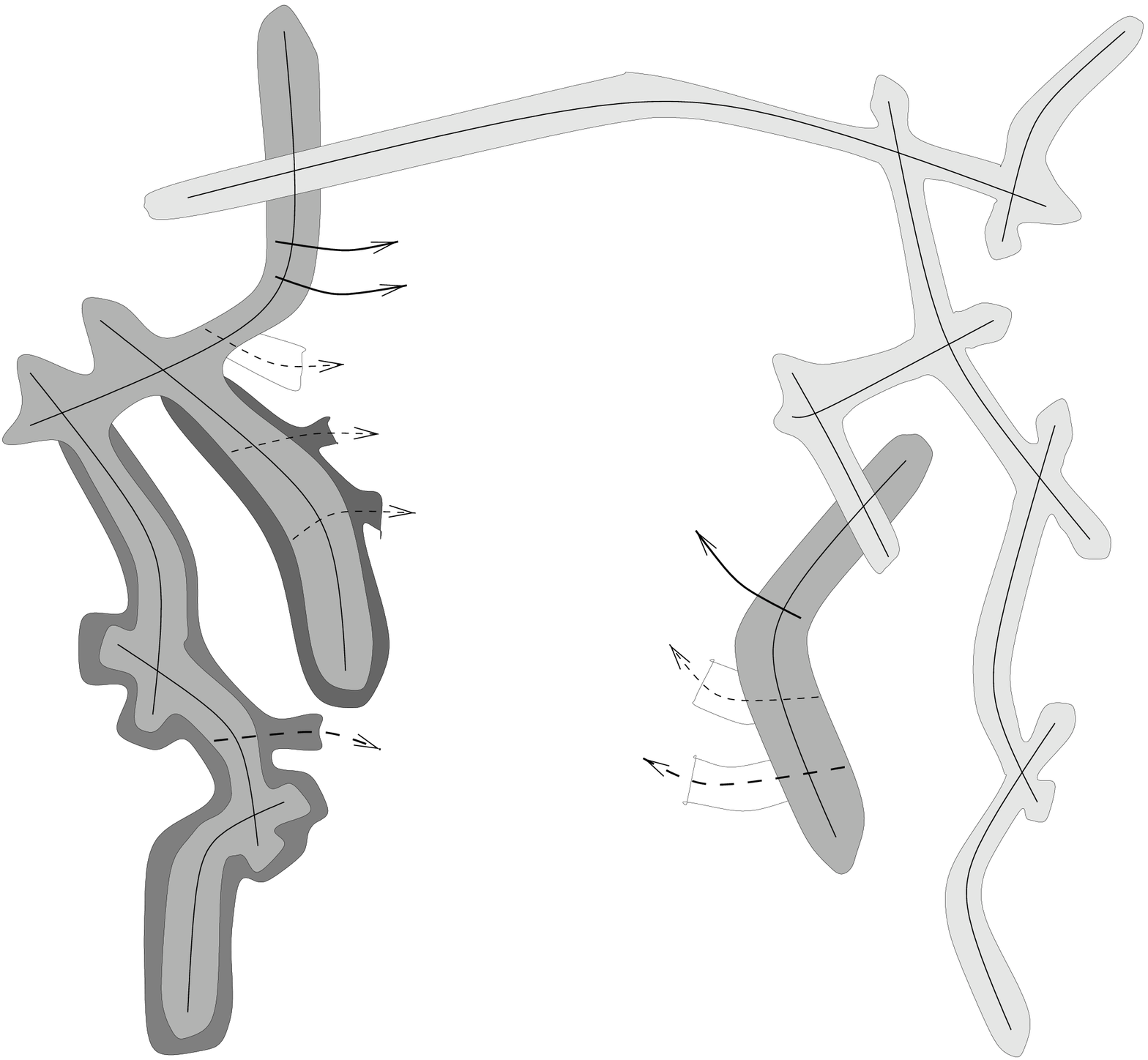}}
\caption{}\label{fig:14}
\end{figure}

Fig.~\ref{fig:14} 
is a schematic picture of the compactification divisor in this case.
We have shaded the domains that are removed in constructing the
manifold $\D$ above:
\begin{itemize}
\item Removing the lightest shaded region removes everything whose
  image under $\fbar$ lies outside the disk $D^2_s(0)$.
\item Removing the next lightest shaded region then removes $|y|>r$,
  to give the 4-disk that we called $D$.
\item Finally, removing the dark regions then removes the components
  of sets  $f^{-1}\bigl(\interior D^2_\epsilon(c)\bigr)
  \cap\bigl\{|y|>r'\bigr\}$ where a fiber $f^{-1}(c)$ is
  irregular at infinity (we have also indicated the components 
  where these fibers are regular at infinity, which are not removed).
\end{itemize}
Thus $\D$ is represented by what has been left white.  Recall that the
parameters in this construction are chosen in the order $s$
(sufficiently large), $r'$ (sufficiently large), $\epsilon$
(sufficiently small), $r$ (sufficiently large). Thus the second
lightest region, which removes a neighborhood of the horizontal and
finite curves of $\Div$, is in fact much the thinnest, although we have
pictured the regions all of comparable size. Most of the boundary
between the dark regions and $\D$ is parallel to fibers of $f$, with
just the small parts near the intersection with the irregular fibers
being transverse to fibers of $f$ (these parts are solid tori).

In \cite{neumann-irreg} it is shown that the splice diagrams at
infinity are derived from the plumbing graph as above, with the
parts of the splice diagram having respectively positive, zero, or
negative total linking weights corresponding to respectively infinite, 
horizontal, or finite curves of $\Div$. We thus see that the tori along
which the splice decomposition of $\partial \D$ occurs are as claimed
in Theorem \ref{th:dim2}, completing the proof (in the picture
these are the places where two different grey-tones meet white).\qed

We close with a comment about the minimal Seifert surface for an
irregular link at infinity.  We first describe how one can see such a
Seifert surface in terms of the construction in the above proof.

Choose a generic line interval $I$ from a point $x$ of $\partial
D^2_\epsilon(c)$ to a point $y$ of $\partial D^2_s(0)$. We can assume,
by choosing $\epsilon$ sufficiently small, that $I$ does not intersect
any of the $\epsilon$-disks around irregular values except at its end
point $x$.  Then $V:=f^{-1}(I)\cap\partial \D$ will be a Seifert
surface for the link at infinity of $f^{-1}(c)$ (considered as a link
in $\partial \D$).  This $V$ is the union of Seifert surfaces for each
of the splice components described in the proof (we work with the
splice decomposition before the elimination of the redundant toral
annulus components):
\begin{itemize}
\item $f^{-1}(y)$ is a fiber of the fundamental multilink and is
  the compact core of a regular fiber.
\item $f^{-1}(x)$ is the disjoint union of the Milnor fibers
  $f^{-1}(x)\cap T_i(c)$ at infinity for $f^{-1}(c)$; these are the
  fibers of fibered splice components corresponding to the $T_i(c)$.
\item $f^{-1}(I)-\Interior(f^{-1}\{x,y\})$ is a union of annuli giving
  Seifert surfaces for the Seifert fibered pieces.
\end{itemize}
Since these are minimal Seifert surfaces for these splice components,
$V$ is a minimal Seifert surface for our link at infinity (see Theorem
3.3 of \cite{eisenbud-neumann}).

Note that the complement of the above minimal Seifert surface $V$ in
the boundary of $f^{-1}(I)\cap \D\cong F\times I$ is the result $F_0$
of removing the Milnor fibers $F_i(c)$ from $f^{-1}(x)\cap \D$. This
$F_0\subset \D$ has boundary isotopic to the link at infinity that we
are considering and it realizes the minimal slice genus of this link,
by the solution of the Thom Conjecture.  The minimal Seifert surface
$V$ is the result of pasting doubles of the Milnor fibers at infinity
onto boundary components of a copy of $F_0$.

Summarizing, our link at infinity, as a link in the boundary of a
4-ball $D^4$, has the property: there is a 2-manifold $F$ containing a
sub-2-manifold $F_0$ and an embedding $F\times [0,1]\subset D^4$ such
that $F_0\times \{0\}\subset D^4$ is a minimal slice surface for the
link while the rest of $\partial (F\times [0,1])$ lies in $\partial
D^4$ and is a minimal Seifert surface for the link.  This is a
presumably already very special property of links at infinity of
affine curves, not shared by general links.  Also special is the fact that
$F\times\{1\}\subset \partial D^4$ and the components of $(F-\interior
F_0)\times \{0\}$ are the fibers of fibered splice components of the
link.


\begin{thebibliography}{99}
  
  \def\title{}
  
\bibitem{artal-cassou-dimca} Artal Bartolo, E., Cassou-Nogu\`{e}s, P.
  and Dimca, A., \title{Sur la topologie des polynomes complexes.}
  Proceedings of Oberwolfach Singularities Conference, (1996),
  Brieskorn Festband, Editors: V.I.Arnold, G.-M.Greuel and J.H.M.
  Steenbrink.

\bibitem{ACL}  E. Artal-Bartolo, P. Cassou-Nogu\'es, I. Luengo Velasco, On
polynomials whose fibers are reducible with no critical points, 
appear. Math. Annalen {\bf299} (1994), 477--490.

\bibitem{broughton1} Broughton, S. A., On the topology of polynomial
  hypersurfaces, Proc. AMS Symp. Pure Math. {\bf40,I} (1983),
  165--178.
  
\bibitem{broughton2} Broughton, S. A., Milnor number and the topology
  of polynomial hypersurfaces, Inv. Math. {\bf92} (1988), 217--241.
  
\bibitem{dimca} Dimca, A., \title{Monodromy at infinity for
    polynomials in two variables}, Preprint, (1998).
  
\bibitem{dimca-nemethi} Dimca, A., Nemethi, Thom Sebastiani
  construction and monodromy of polynomials, Pr\'epublication {\bf 98}
  (1999), Laboratoire de Math. Pures de Bordeaux C.N.R.S.
  
\bibitem{durfee-var} Durfee, A., Fibered knots and algebraic
  singularities, Topology {\bf13} (1974), 47--59.

\bibitem{durfee} 
 Durfee, Alan H. Five definitions of critical point at
 infinity. Singularities (Oberwolfach, 1996), Progr. Math.,
 {\bf 162}, 345--360

  
\bibitem{eisenbud-neumann} Eisenbud, D. and Neumann, W.D.,
  \emph{Three-dimensional link theory and invariants of plane curve
    singularities.}  Ann. Math. Stud. {\bf 110}, Princeton.  Princeton
  Univ. Press (1985).
  
\bibitem{hirsch} Hirsch, M. \emph{Differential Topology}, Graduate
  Texts in Math. {\bf33} (Springer Verlag, 1976).
  
\bibitem{jankins-neumann} M. Jankins and W.D. Neumann, Lectures on
  Seifert manifolds, Brandeis Course Notes {\bf 2} (1983).

\bibitem{lamotke} Lamotke, Klaus, Die Homologie isolierter
  Singularit\"aten. Math. Z. {\bf143} (1975), 27--44.

\bibitem{milnor} Milnor, J., \emph{Singular points of complex
    hypersurfaces}, Ann. Math. Stud. {\bf 101}, Princeton University
  Press, (1968).
  
\bibitem{nemethi-zaharia} N\'emethi, A., Zaharia, A., On the
  bifurcation set of a polynomial and Newton boundary, Publ. RIMS
  {\bf26} (1990), 681--689.
  
\bibitem{neumann-splice} Neumann, W.D., Splicing algebraic links, in
  \emph{Complex Analytic Singularities}, Advanced Studies in Pure
  Math. {\bf8} (1986), 349--361.

\bibitem{neumann-inv} Neumann, W.D., \title{Complex algebraic curves
    via their links at infinity}, Invent. Math. {\bf 3}, (1989),
  445-489.
  
\bibitem{neumann-irreg} Neumann, W.D., \title{Irregular links at
    infinity of complex affine plane curves}, Quarterly
  J. Math. {\bf{50}} (1999), 301--320.
  
\bibitem{neumann-norbury} Neumann, W.D. and Norbury, P.,
  \title{Monodromy and vanishing cycles of complex polynomials}, Duke
  Math. J. (to appear).

  
\bibitem{neumann-rudolph} Neumann, W. and Rudolph, L.
  \title{Unfoldings in knot theory}, Math. Ann. {\bf 278}, (1987),
  409-439.  \emph{Corrigendum}, Math. Ann. {\bf 282}, (1988), 349-351.
  
\bibitem{parusinski} Parusi\'nski, A. On the bifurcation set of a
  complex polynomial with isolated singularities at infinity,
  Compositio Math. {\bf97} (1995), 369--384.

\bibitem{parusinski2} Parusi\'nski, A. A note on singularities at
  infinity of complex polynomials, \emph{Simplectic singularities and
  geometry of gauge fields}, Banach Center Publ. {\bf39} (1997),
  131--141.
  
\bibitem{pham} Pham, F., Vanishing homologies and the $n$ variable
  saddlepoint method, AMS Proc. Sympos. Pure Math. {\bf40,II} (1983),
  319--333.

\bibitem{siersma-tibar} Siersma, D., Tib\v ar, M. Singularities at
  infinity and their vanishing cycles, Duke Math. J. {\bf80} (1995),
  771--783.
  
\bibitem{suzuki} Suzuki, M.  \title{Propri\'{e}t\'{e}s topologiques
    des polyn\^{o}mes de deux variables complexes et automorphismes
    alg\'{e}briques de l'espace $C^2$}, J. Math. Soc. Japan {\bf 26},
  (1974), 241-257.
  
\bibitem{suzuki77} Masakazu Suzuki, Sur les op\' erations holomorphes
  du groupe additif complexe sur l'espace de deux variables complexes.
  (French) Ann. Sci. \' Ecole Norm. Sup. (4) {\bf10} (1977), 
  517--546.
  
\bibitem{tibar} Tib\v ar, M., On the monodromy fibration of
  polynomial functions with singularities at infinity,
  C. R. Acad. Sci. Paris, {\bf324}, S\'erie I (1997), 1031--1035.
  
\bibitem{tibar-Wallpaper} Tib\v ar, M., Regularity at infinity of real
  and complex polynomial functions, \emph{Singularity Theory, C.T.C.
    Wall Anniversary Volume}, (Cambridge U. Press).

\end{thebibliography}
\end{document}